\documentclass[preprint,12pt]{elsarticle}

\usepackage{graphicx}
\usepackage{amssymb}

\usepackage{lineno}

\usepackage{graphicx}
\usepackage[utf8]{inputenc}
\usepackage{lmodern}
\usepackage[T1]{fontenc}
\usepackage{lscape}
\usepackage{amsmath}	
\usepackage{amssymb}
\usepackage{cases}
\usepackage{graphicx}
\usepackage{caption}
\usepackage{comment}

\usepackage{makecell}
\usepackage{multirow}
\usepackage{amsfonts}
\usepackage{float} 
\usepackage{diagbox}
\usepackage{fullpage}
\usepackage{url}
\usepackage{rotating}
\usepackage{eurosym}
\usepackage{wrapfig}
\usepackage[final]{pdfpages} 
\usepackage{epstopdf} 
\usepackage[a4paper]{geometry}
\usepackage[nottoc]{tocbibind}
\usepackage{placeins}
\usepackage{float}
\usepackage{listings}
\usepackage{color, colortbl}
\usepackage{hyperref}
\usepackage{xcolor}
\usepackage{graphicx, caption, subcaption}

\usepackage{sidecap}
\hypersetup{
colorlinks,
citecolor=black,
filecolor=black,
linkcolor=black,
urlcolor=black
} 

\newcommand{\argmin}{\operatornamewithlimits{argmin}}

\newcommand{\bx}{\textbf{x}}

\newcommand{\by}{\textbf{y}}
\newcommand{\bB}{\textbf{B}}
\newcommand{\bH}{\textbf{H}}
\newcommand{\bR}{\textbf{R}}
\newcommand{\bA}{\textbf{A}}

\newcommand{\bD}{\textbf{D}}
\newcommand{\bM}{\textbf{M}}

\newcommand{\bY}{\textbf{Y}}
\newcommand{\bC}{\textbf{C}}
\newcommand{\bL}{\textbf{L}}

\newcommand{\bS}{\textbf{S}}
\newcommand{\bI}{\textbf{I}}

\journal{Journal of computational science}

\begin{document}

\begin{frontmatter}


\title{Observation data compression for variational assimilation of dynamical systems}



\author{Sibo Cheng$^{1,2,3}$, Didier Lucor$^{2}$, Jean-Philippe Argaud$^{3}$
}

\address{   \small $^{1}$ Data Science Institute, Department of Computing, Imperial College London, UK\\
        \small $^{2}$ Université Paris-Saclay, CNRS, Laboratoire Interdisciplinaire des Sciences du Numérique, France\\
        \small $^{3}$ EDF R\&D Saclay, France\\}

\begin{abstract}
 Accurate estimation of error covariances (both background and observation) is crucial for efficient observation compression approaches in data assimilation of large-scale dynamical problems. We propose a new combination of a covariance tuning algorithm with existing PCA-type data compression approaches, either observation- or information-based, with the aim of reducing the computational cost of real-time updating at each assimilation step. Relying on a local assumption of flow-independent error covariances, dynamical assimilation residuals are used to adjust the covariance in each assimilation window. The estimated covariances then contribute to better specify the principal components of either the observation dynamics or the state-observation sensitivity. The proposed approaches are first validated on a shallow water twin experiment with correlated and non-homogeneous observation error. Proper selection of flow-independent assimilation windows, together with sampling density for background error estimation, and sensitivity of the approaches to the observations error covariance knowledge, are also discussed and illustrated with various numerical tests and results. The method is then applied to a more challenging industrial hydrological model with real-world data and non-linear transformation operator provided by an operational precipitation-flow simulation software.
\end{abstract}


\end{frontmatter}


\section{Introduction}
Data assimilation (DA) is applied in a wide range of industrial problems, such as numerical weather prediction (NWP) \cite{Rabier2005}, hydrology, fire forecasting \cite{rochoux2014towards} or nuclear engineering \cite{Gong2020}.  Recently, DA methods have also been used to COVID-19 pandemic analysis, including predicting disease diffusion and proposing optimal vaccination strategies (\cite{cheng2021}). DA algorithms are often used in dynamical systems for continuously updating state estimation/prediction. They have recently made their way to other fields such as biomedical applications \cite{lucor_esaim2018} or quantitative economics \cite{Nadler2019}. These methods rely on a weighted combination of different sources of noisy information, including prior numerical estimation (also known as background states) and real-time observations, to improve field reconstruction or parameters calibration. DA methods are often used to deal with problems of large dimensions, especially in NWP \cite{collard2010}, \cite{Fowler2019} (up to $10^9$) or in geoscience \cite{Carrassi2017}, leading to computational difficulty for real-time updating, if not infeasible. Several strategies for optimizing the computational cost have been developed, including graph-based domain localization \cite{cheng2020}, observations selection \cite{Cioaca2014}, matrix decomposition \cite{Nino-Ruiz2019} or reduced-order Kalman Filter \cite{Hoteit2002}. It is also a common practice to combine DA algorithms with classical dynamical system reduction techniques, such as the Proper Orthogonal Decomposition (POD) or the Empirical Interpolation Method  (EIM e.g. \cite{helin2018}). Most of these methods rely on either precise knowledge of state variables (e.g. modes in POD) or strong prior assumptions (e.g. cut-off radius in domain localization \cite{Waller2017}). Meanwhile, with the increase of available observation precision in DA applications, the observation data compression via low-rank approximation methods has been continuously studied for alleviating the computational cost, especially in a sequential data assimilation chain.  These methods, which consist of extracting principal information in observation data, have been widely applied in various branches of engineering, especially for high dimensional problems. An important advantage of observation compression, regarding other methods that directly reduce the state space dimension, is that no extra operation/knowledge of the state dynamics is required, making the compression error more controllable and estimable. Two classical compression methods are discussed and implemented in this work: the POD-type projection by extracting principal components in the observation dynamic \cite{Matricardi2014} and the information-based compression based on the information entropy analysis \cite{Fowler2019}. The latter aims to select the most impacting observations to the analyzed state by calculating the prior-posterior information entropy gap. Since the noises are introduced by prior errors in DA systems, the information entropy estimation relies on both background and observation error covariance matrices. \\

For both observation- and information-based approaches, the data compression is carried out with a noise-normalized dataset \cite{collard2010}, \cite{Fowler2019}. The knowledge of prior error covariances thus becomes crucial for applying these methods. However, the specification of these covariances, especially the background matrix, remains one of the most challenging problems in data assimilation due to the high dimension of the problem and limited prior data \cite{Fisher2003},\cite{cheng2019}. Much attention was given to improving the error covariance specification in dynamical data assimilation models, particularly by the meteorological society. Several methods have been developed to this end, such as the NMC approach \cite{F.Parrish1992}, the DI01 \cite{Desroziers01} iterative method and the Desroziers estimation \cite{Desroziers2005}.  In this paper, we focus on the latest. Unlike some other methods (e.g. \cite{Desroziers01}, \cite{cheng2019}), the Desroziers estimation does not depend on the specific structure of the error covariances, and it provides a non-parametric estimation of full covariances as output of the algorithm. Based on the residual analysis in variational assimilation, this approach has been widely applied in industrial problems, especially in NWP. Recent works of \cite{Bathmann2018} prove its convergence in the ideal case. Another considerable strength of the Desroziers estimation is that dynamic residual data can be used for the covariance estimation.  For this reason, a huge ensemble size is not required for high dimensional problems, unlike, for instance, in the NMC method. \\

In this paper, based on the Desroziers estimation, we have introduced the concept of piecewise covariance estimation for both observation- and information-based compression strategies. We apply the Desroziers method to estimate error covariances in a fixed time range, also known as the flow-independent window where the error covariances are supposed to be time-invariant.  Therefore, the choice of the flow-independent window and the residual samplings play an essential role in this algorithm. The window size should be sufficiently long to gather enough time-variant sampling but not too long to consider the error covariances, especially the background matrix, being constant. \\

The observation- and information-based (with piecewise covariances estimation) data compression are first implemented in a twin experiment framework using 2D shallow water equations with a linear transformation operator. The observation covariance is supposed to be perfectly known \textit{a priori}. The two approaches with different choices of flow-independent windows are compared in this model while changing the truncation parameter. Numerical results show that the observation-based (POD-type) compression is in general over-performed by the information-based approach and that a non-balanced sampling in piecewise covariance estimation results in a less optimal compression. We then apply these methods to a real-world hydrological model to improve river flow prediction/reanalysis by correcting historical daily precipitation measures \cite{cheng2020b}. Both the precipitation and the river flow data are spatially distributed. The physical simulation is performed using the operating MORDOR-TS software \cite{Rouhier2017}, developed by EDF and the study area is around the Tarn river, in the south of France. The precipitation-flow simulation is carried out through conceptual
watersheds modeling, which ensures its high computational efficiency. In this hydrological application, both the background and the observation matrices are estimated using the Desroziers method with daily observed flow data for around 10 years (1990 to 2000). Results show that in this industrial application where both $\bB$ and $\bR$ are not well known, the performance of the information-based strategy is similar to the one of observation-based.\\

The paper is organized as follows. In section \ref{sec:DA}, the principle and the notation of data assimilation are briefly introduced. We then introduce the observation- and information-based compression strategies in section \ref{sec: metho}. The applications of 2D shallow water twin experiments and an industrial hydrological model are shown respectively in section \ref{sec:sw} and \ref{sec:hydro}. We finish the paper with a discussion.

\section{Variational data assimilation}
\label{sec:DA}
The objective of data assimilation algorithms is to improve the estimation of some physical fields or parameters $\textbf{x}$ based on two sources of information: a prior simulation/forecast $\textbf{x}_b$ and an observation vector $\textbf{y}$. The theoretical value of the current state is denoted by a vector $\textbf{x}_\textrm{true}$, also known as the true state. 
Variational DA algorithms aim  to find an optimally weighted compromise between the prior estimation $\textbf{x}_b$ and the observation $\textbf{y}$ by minimising the cost function $J$ defined as
\begin{align}
    J(\textbf{x})&=\frac{1}{2}(\textbf{x}-\textbf{x}_b)^T\textbf{B}^{-1}(\textbf{x}-\textbf{x}_b) + \frac{1}{2}(\textbf{y}-\mathcal{H}(\textbf{x}))^T \textbf{R}^{-1} (\textbf{y}-\mathcal{H}(\textbf{x})) \label{eq_3dvar}\\
   &=\frac{1}{2}||\textbf{x}-\textbf{x}_b||^2_{\textbf{B}^{-1}}+\frac{1}{2}||\textbf{y}-\mathcal{H}(\textbf{x})||^2_{\textbf{R}^{-1}}
\end{align}
 where  $\mathcal{H}$ denotes the transformation operator from the state space to one of the observations. $\textbf{B}$ and $\textbf{R}$ are the associated error covariance matrices, i.e.
 \begin{align}
     \textbf{B} = \textrm{Cov}(\epsilon_b, \epsilon_b), \quad
     \textbf{R} = \textrm{Cov}(\epsilon_y, \epsilon_y),
 \end{align}
 where
  \begin{align}
     \epsilon_b = \textbf{x}_b - \textbf{x}_\textrm{true}, \quad
     \epsilon_y = \mathcal{H}(\textbf{x}_\textrm{true})-\textbf{y}.
 \end{align}
Thus the inverse of these covariance matrices (i.e. $\textbf{B}^{-1}, \textbf{R}^{-1}$) represents the weights of these two information sources in the objective function. Prior errors $\epsilon_b, \epsilon_y$ are supposed to be centered Gaussian, characterised by the error covariance matrices, i.e.
 \begin{align}
     \epsilon_b \sim \mathcal{N} (0, \textbf{B}), \quad
     \epsilon_y \sim \mathcal{N} (0, \textbf{R}).
 \end{align}
 
 The optimization problem of Eq.~\ref{eq_3dvar}, so called three-dimensional variational (3D-Var) formulation, is a general representation of variational assimilation while the model error is not considered. The output of Eq.~\ref{eq_3dvar} is denoted as $\bx_a$, i.e.
  \begin{align}
    \bx_a = \underset{\bx}{\argmin} \Big(J(\textbf{x})\Big). \label{eq:argmin}
 \end{align}
 If $\mathcal{H}$ can be approximated by some linear operator $\bH$,  Eq.~\ref{eq:argmin} can be solved via BLUE (Best Linearized Unbiased Estimator) formulation,
   \begin{align}
    \textbf{x}_a &= \textbf{x}_b+\textbf{K}(\textbf{y}-\textbf{H} \textbf{x}_b) \\
    \bA &= (\textbf{I}-\textbf{K}\textbf{H})\textbf{B} \label{eq:BLUE}
 \end{align}
 where $\textbf{A} = \textrm{Cov}(\textbf{x}_a-\textbf{x}_\textrm{true})$ is the analyzed error covariance and the $\textbf{K}$ matrix, given by
 \begin{equation}
	\textbf{K}=\textbf{B} \textbf{H}^T (\textbf{H} \textbf{B} \textbf{H}^T+\textbf{R})^{-1}
	\label{eq:Kgain_BLUE}
\end{equation}
is so called the Kalman gain matrix.
In the rest of this paper, we denote $\bH$ as the linearized transformation operator. The case when $\mathcal{H}$ is non-linear is more challenging for finding the minimum of Eq.~\ref{eq_3dvar}, especially for high-dimensional problems. The resolution involves often gradient descent algorithms (relying on algorithms such as "L-BFGS-B" \cite{Fulton2000} and on adjoint-based \cite{Cioaca2014} numerical techniques. \\

Variational assimilation algorithms could be applied to dynamical systems through sequential applications using a transition operator $\mathcal{M}_{t^k \rightarrow t^{k+1}}$ (from time $t^k$ to $t^{k+1}$), where
\begin{align}
    \textbf{x}_{t^{k+1}} = \mathcal{M}_{t^k \rightarrow t^{k+1}} (\textbf{x}_{t^k }).
\end{align}
The forecasting thus depends on the knowledge of transition operator  $\mathcal{M}_{t^k \rightarrow t^{k+1}}$ and the corrected state at the current time $ \textbf{x}_{a, t^k }$. Typically, the current background state is often given by the forecasting from the previous step, i.e.
\begin{align}
    \textbf{x}_{b,t^k} = \mathcal{M}_{t^{k-1} \rightarrow t^{k}} (\textbf{x}_{a,t^{k-1} }).
\end{align}
Obviously, a more accurate reanalysis $\textbf{x}_{a,t^{k-1}}$ leads to a more reliable forecasting $\textbf{x}_{b,t^k}$. It is known that as long as the transformation operator $\mathcal{H}$ and the transition operator $\mathcal{M}$ are linear, the analysis based on the variational method and the Kalman filter results in the same forecasting \cite{Carrassi2017}, for dynamical (4D-Var) assimilation problems. Theoretically, the evolution of the $\bB$ matrix could also be estimated thanks to the transition operator. However, in practice, the pefect knowledge of $\mathcal{M} $ is often unavailable. Much attention is given to quantify the model error in assimilation, for example, in weak-constraint \textit{4D-VAR} \cite{UBOLDI2000}. Recent work of \cite{Brajard2020} involves deep learning techniques to improve the estimation of  $\mathcal{M}_{t^{k-1} \rightarrow t^{k}}$.

\section{Observation data compression}
\label{sec: metho}
DA algorithms are often used to perform real-time corrections of dynamical systems with large dimensions, leading to an essential requirement of computational efficiency. In this work, we are interested in a low-rank approximations of the observation vector which can reduce the cost of real-time updating in DA algorithms.
\subsection{Observation-based compression (OC)}
 The works of \cite{Antonelli2004} and \cite{Tobin2006} are based on a PCA-type reduction of the observation dynamics. More precisely, a set of $n_{obs}$ observation snapshots is represented by a matrix $\bY \in \mathbb{R}_{[\textrm{dim}(\by) \times n_{obs}]}$ where each column $\bY [:,.]$ represents an individual observation vector of dimension $m$ at a fixed time $t_i$, i.e.
\begin{align}
    \bY [:,i] = \by_{t=t_i}.
\end{align}
Thus $\bY$ describes the evolution of the observation vector $\by$ including observation error. We work with the error-normalized data $\bR^{-1/2} \bY$ \cite{collard2010} whose empirical covariance $\bC$ can be written and decomposed as

\begin{align}
    \bC = \frac{1}{n_{obs}-1} \bR^{-1/2} \bY \bY^T \bR^{-1/2} = \tilde{\bL} \tilde{\bD} \tilde{\bL}^T
\end{align}
where the columns of $\tilde{\bL}$ are the principal components and $\tilde{\bD}$ represents the associated eigenvalues in a decreasing order. This decomposition is known as the principal component analysis (PCA) decomposition. We can construct a projection operator $\tilde{\bL}_q$ with minimum loss of information (represented by eigenvalues in the covariance matrix) by simply keeping the $q$ first columns in $\tilde{\bL}$. $q$ is also known as the truncation parameter. In fact, this projection operator can also be obtained by a singular value decomposition (SVD), without computing the full covariance matrix $\bC$, i.e.
\begin{align}
    \bR^{-1/2}  \bY = \tilde{\bL}_q \tilde{\Sigma} \tilde{\textbf{V}_q}^T \label{eq:SVD_Y}
\end{align}
where $\tilde{\textbf{L}_q}$ and $\tilde{\textbf{V}_q}$ are orthogonal matrices, i.e. $\tilde{\textbf{L}_q}^T \tilde{\textbf{L}_q} = \tilde{\textbf{V}_q}^T \tilde{\textbf{V}_q} = \bI $ and $\tilde{\Sigma}\tilde{\Sigma}^T = \tilde{\bD}$ since all eigenvalues are non negative. The assumption is made for the observation error covariances to be constant (flow-independent), which is a common pratice in data assimilation (e.g \cite{collard2010}). For each DA optimization, instead of updating with the full observation vector $\by$, the correction is made with the reduced observation
\begin{align}
    \tilde{\by}_q =  \tilde{\bL}^T_q \bR^{-1/2} \by. \label{eq: reconstruction}
\end{align}

The new observation error covariance $\tilde{\bR}$ and the new state-observation transformation operator $\tilde{\mathcal{H}}$ can be written as
\begin{align}
    \tilde{\bR}_q =  \tilde{\bL}^T_q \bR^{-1/2} \bR \bR^{-1/2} \tilde{\bL}_q = \bI_q, \quad \tilde{\mathcal{H}}_q = \tilde{\bL}^T_q \bR^{-1/2} \circ \mathcal{H}. \label{eq:OC_R}
\end{align}

The DA algorithm can then be performed on $(\bx_b, \tilde{\by}_q, \bB, \tilde{\bR}_q, \tilde{\mathcal{H}}_q)$ instead of $(\bx_b, \by, \bB, \bR, \mathcal{H})$. This method could be seen as a classical POD approach applied to error-normalised observation data by extracting modes of higher variances against time. It is pointed out by \cite{Antonelli2004} and \cite{collard2010} that performing PCA on noise-normalised observation data  can improve the method efficiency and reduce the impact of observation error during the compression procedure.

\subsection{Information-based compression (IC)}
The observation-based data reduction retains the principal directions of the observation dynamic. However, these directions are not necessarily the most impacting in state correction. A continuous effort has been devoted to quantify and compute the sensitivity of the analysis states to the observations (e.g. \cite{Cioaca2014}), which leads to a more refined observation compression in DA. More precisely, this sensitivity may be expressed by the influence matrix $\bS$ \cite{Cardinali2204}, defined as
\begin{align}
    \bS = \frac{\partial \mathcal{H}(\bx_{a})}{\partial \bx_{a}} = \textbf{K}^T \bH^T.
\end{align}
 According to \cite{Fowler2019}, the information given by the influence matrix can be roughly quantified via two indicators, the degree of freedom for signal (DFS) which represents the prior-posterior mutual information and the entropy reduction (ER) which represents the evolution of Shannon information content, respectively defined as
\begin{align}
    \textrm{DFS} &= \mathbb{E} [(\bx_a - \bx_b)^T \bB^{-1} (\bx_a - \bx_b)] = Tr (\bS) \label{eq:info2}\\
    \textrm{ER} & = H(\bx) - H(\bx|\by) = -\frac{1}{2} \ln \big(\textit{det} (\bI - \bS)\big)\label{eq:info}
\end{align}
where $H$ is the entropy of a distribution, noted here $H(x)$ for simplicity. Eqs.~\ref{eq:info2} and \ref{eq:info} are derived for a centred Gaussian vector $\mathbf{x}$. For both measures, we observe that observations associated with the largest eigenvalues of $\bS$ have the greatest information content. Using an intermediate matrix $\bM = \bR^{-1/2} \bH \bB ^{1/2}$, Eqs.~\ref{eq:info2}-\ref{eq:info} could be rewritten as
\begin{align}
    \textrm{DFS} &= \textrm{Tr}(\bM \bM^T (\textbf{I}+\bM \bM^T)^{-1}) \\
    \textrm{ER} &= \frac{1}{2} \ln \big(det (\textbf{I}+\bM \bM^T)\big). \label{eq: M}
\end{align}
  
As stated in the work of \cite{Migliorini2013}, the observation projection operator which minimizes the information loss is given by $\hat{\bL}_q \bR^{-1/2}$, where $\hat{\bL}_q$ is the matrix whose columns contain the eigenvectors of $\bM \bM^T = \bR^{-1/2} \bH \bB \bH^T \bR^{-1/2}$. DA algorithms could then be performed with 
\begin{align}
  \hat{\by}_q =  \hat{\bL}^T_q \bR^{-1/2} \by, \quad
    \hat{\bR}_q =  \bI_q, \quad \hat{\mathcal{H}}_q = \hat{\bL_q}^T_q \bR^{-1/2} \circ \mathcal{H}. \label{eq:IC_R}
\end{align}
We remind that, from the computational point of view, the only difference between the OC and IC is the way the low-rank projection $\bL_q$ is obtained. For both approaches, the specification of error covariance matrices (either background or observation) is crucial to provide an efficient compression. On the other hand, data compression strategies can reduce the computational cost of covariance  tuning methods, especially for multidimensional and multivariate problems. Therefore, the precise knowledge of $\textbf{H} \textbf{B} \textbf{H}^T$ and $\bR$ is crucial for this method. However, as pointed out by \cite{Fowler2019}, the condition number of the analysis covariance matrix $\bA$ can be higher when using IC approach compared to performing DA with the full observation data set. Therefore, the risk of matrix ill-conditioning is worth monitoring when applying this compression method.

\subsection{Optimal truncation parameter for compression methods}
\label{sec:stop}
The determination of the truncated parameter $q$, i.e. number of modes kept in the reduced space, is crucial in data compression. The choice of the threshold often depends on available data \cite{Cangelosi2007}.  Several criteria were considered, such as the information losing rate $E_q$ and the matrix conditioning \textit{a posteriori} $\mu_q$, defined as
\begin{align}
 E_q &= \frac{||\Sigma-\phi_q||_\infty}{||\Sigma||_\infty} = 1 -\frac{\sigma_{q-1}}{\sigma_q} \label{indc1}\\
 \mu_q &= \frac{\sigma_1}{\sigma_q} \label{indc2}
\end{align}
where $\Sigma$ is the diagonal matrices with all eigenvalues of the covariance matrix and $\sigma_{i, i=1..}$ represent the associated real eigenvalues in the decreasing order of absolute value. 
According to the study of \cite{Arcucci2018}, an optimal choice of the truncation parameter can be obtained by combining the two previous indicators, with an objective function $f$, defined as
\begin{align}
    f(\sigma_q) &= E_q + \mu_q  = \frac{\sigma_{q} \sigma_{q-1} + \sigma_1}{\sigma_1 \sigma_{q}}. \label{eq:q_op}
\end{align}
Assuming $||\sigma_{q}-\sigma_{1}||>>||\sigma_{q}-\sigma_{q-1}|| $, one could easily prove that Eq.~\ref{eq:q_op} achieves the minimum when $\sigma_{q} = \sqrt{\sigma_{1}}$. With this choice, we manage to both reduce the matrix ill-conditioning and remove less significant modes, as proved in real-world DA application\cite{Arcucci2018}. Another advantage of this criteria is that the computation of the full spectrum of covariances is not required. By applying Lancozs-type methods \cite{Lanczos1950}, we can stop the algorithm when the current eigenvalue is inferior to $\sqrt{\sigma_{1}}$.

\section{Piecewise estimation of error covariances}

The $\bR$ matrix is required for both OC and IC approaches. Furthermore, the construction of $\hat{\bL}^T_q$ in IC requires a precise knowledge of the matrix production $\bH \bB \bH^T$. However, the knowledge of both matrices often remains challenging in data assimilation \cite{Fisher2003}. Continuous effort was devoted to improve the error covariances specification \cite{Tandeo2018}, \cite{cheng2019}.  A classical approach based on residual analysis, and later a more complete version are respectively given by \cite{Hollingsworth1989b} and \cite{Desroziers2005}. They show that under the assumption of flow-independent error covariance, i.e. $\bB$ and $\bR$ being invariant against time in a certain period, the following equations hold
\begin{align}
         \bR & =\mathbb{E}\Big[\big(\by-\mathcal{H}(\bx_a)\big)\big(\textbf{y}-\mathcal{H}(\bx_b)\big)^T\Big] \label{eq:d05R}\\
    \bH \bB \bH^T & = \mathbb{E}\Big[\big(\by-\mathcal{H}(\bx_b)\big)\big(\by-\mathcal{H}(\bx_b)\big)^T\Big] - \bR.  \label{eq:hollings} 
\end{align}
Under these hypothesis, combining Eq.~\ref{eq:d05R} and \ref{eq:hollings} leads to 
\begin{align}
          \bH \bB \bH^T & =\mathbb{E}\Big[\big(\mathcal{H}(\bx_a)-\mathcal{H}(\bx_b)\big)\big(\textbf{y}-\mathcal{H}(\bx_b)\big)^T\Big]. \label{eq:d05}
\end{align}
In order to somewhat alleviate these strong hypotheses, a simple idea is to take the expectation operators in Eq.~\ref{eq:d05R}, \ref{eq:hollings} and \ref{eq:d05} in assimilation windows where the flow-independent assumption stands, resulting in a piecewise estimation of both $\bB$ and $\bR$. More precisely, a sequence of estimated background matrices $\bB_{T_i}$ could be computed via residual covariances, where $T_i$ refer to flow-independent periods of $\bB$ in a dynamical system. In other words, $\bB$ is considered as invariant between $t=T_i$ and $t=T_{i+1}$. The estimation of $\bR_{T_i}$, if required, follows the same principle using Eq.~\ref{eq:d05R}. When the knowledge of $\bR$ matrix is precise \textit{a priori}, the estimation of Eq.~\ref{eq:hollings} is privileged because of its lower computational cost since no evaluation of the analyzed state $\bx_a$ is required.  According to \cite{Hollingsworth1989b}, when the observation error is dominated by background error (i.e. $\textrm{Tr}(\bR)<<\textrm{Tr}(\bB)$), $\bH \bB \bH^T$ can be estimated directly by $\mathbb{E}\Big[\big(\by-\mathcal{H}(\bx_b)\big)\big(\by-\mathcal{H}(\bx_b)\big)^T\Big]$. 

By definition, 
\begin{align}
    \bH \bB \bH^T  = \mathbb{E}\Big[\big(\mathcal{H}(\bx_\textrm{true})-\mathcal{H}(\bx_b)\big)\big(\mathcal{H}(\bx_\textrm{true})-\mathcal{H}(\bx_b)\big)^T\Big] \label{eq:HBH}
\end{align}
represents the background error covariances projected in the observation space.Therefore the information-based observation compression which is based on a PCA-type analysis, can also be interpreted as a projection of $\by$ along the directions where the background errors are most important. Recently, it was also reported in the literature (e.g. \cite{Waller2016b}) that the convergence (towards the exact observation matrix) of the iterative method can still be ensured when the background and observation error correlation length-scales are similar, which was contrary to what was previously thought \cite{Chapnik2006}. Although this innovation-based covariance estimation approach has been widely applied in DA applications, some drawbacks have also been noticed. For example, the application of this method in real problems often requires post-processing of the $\bR$ matrix. It is shown in \cite{cheng2020b} that the regularized matrix may converge to some other solution rather than the exact observation matrix.

\section{Shallow water twin experiments}
\label{sec:sw}
\subsection{Experiments set up}
\label{sec:Experiments set up}
For evaluating the performance of different data compression approaches, we set up a twin experiment framework with  a simplified 2D shallow water dynamical model which is frequently used for testing data assimilation algorithms (
 (e.g \cite{Cioaca2014}, \cite{cheng2019}). A cylinder of water is positioned in the middle of the study field of size $20mm \times 20mm$ and released at the initial time $t=0s$ (i.e. with no initial speed), leading to a non-linear wave-propagation. The dynamics of the water level $h$ (in $mm$), as well as horizontal and vertical velocity (in $0.1m/s$) field (respectively denoted as $u$ and $v$), is given by the non-conservative shallow water equations

\begin{align}
    \frac{\partial u}{\partial t}&=-g\frac{\partial}{\partial x}(h)-bu \label{eq: sw} \\
    \frac{\partial v}{\partial t}&=-g\frac{\partial}{\partial y}(h)-bv  \notag\\ 
    \frac{\partial h}{\partial t}&=-\frac{\partial}{ \partial x}(uh)-\frac{\partial}{\partial y}(v h)  \notag \\
    u_{t=0} &= 0 \notag \\
    v_{t=0} &= 0 \notag
\end{align}

\begin{figure}
  \centering
  \subfloat[]{\includegraphics[width=2.2 in]{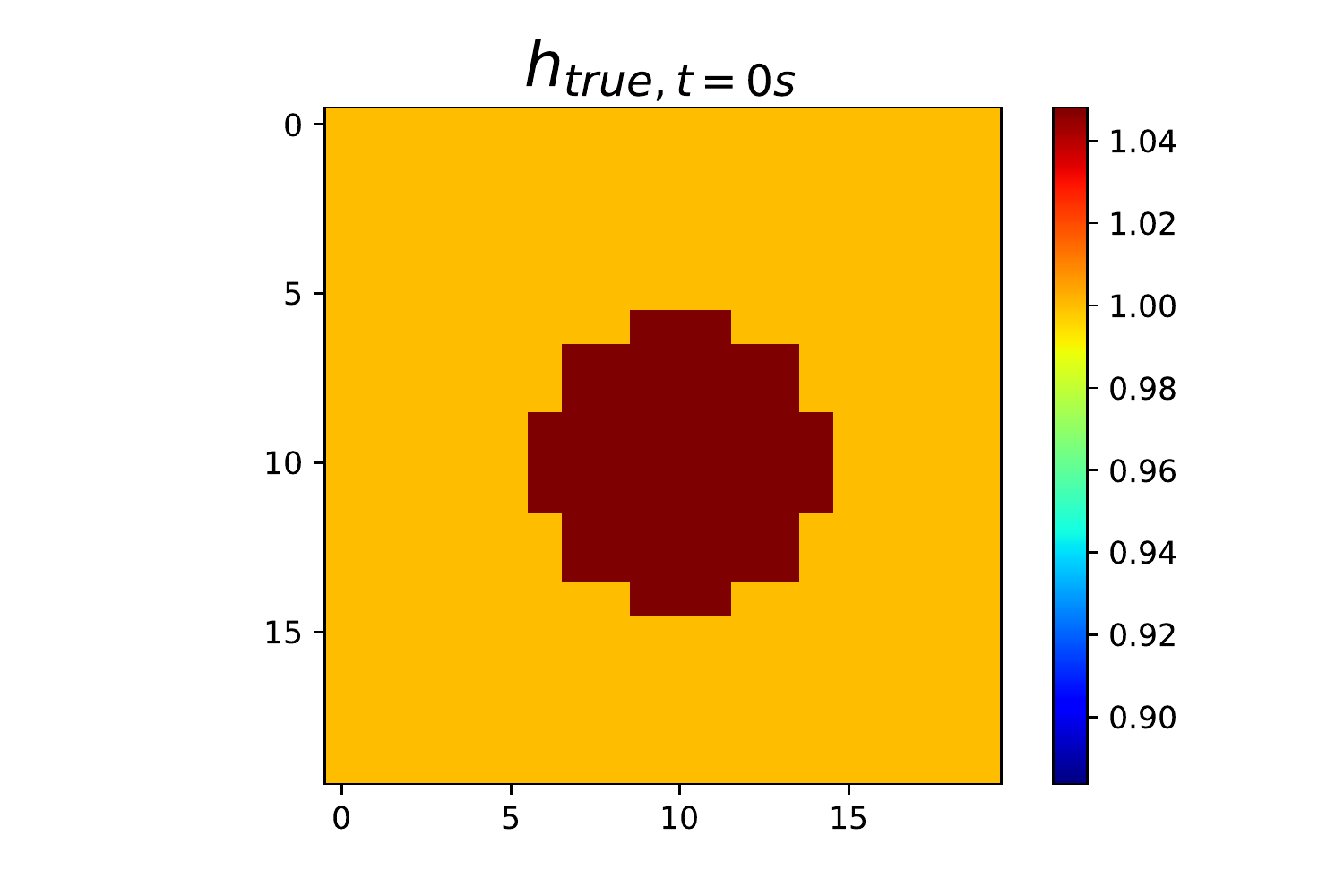}}
   \subfloat[]{\includegraphics[width=2.2 in]{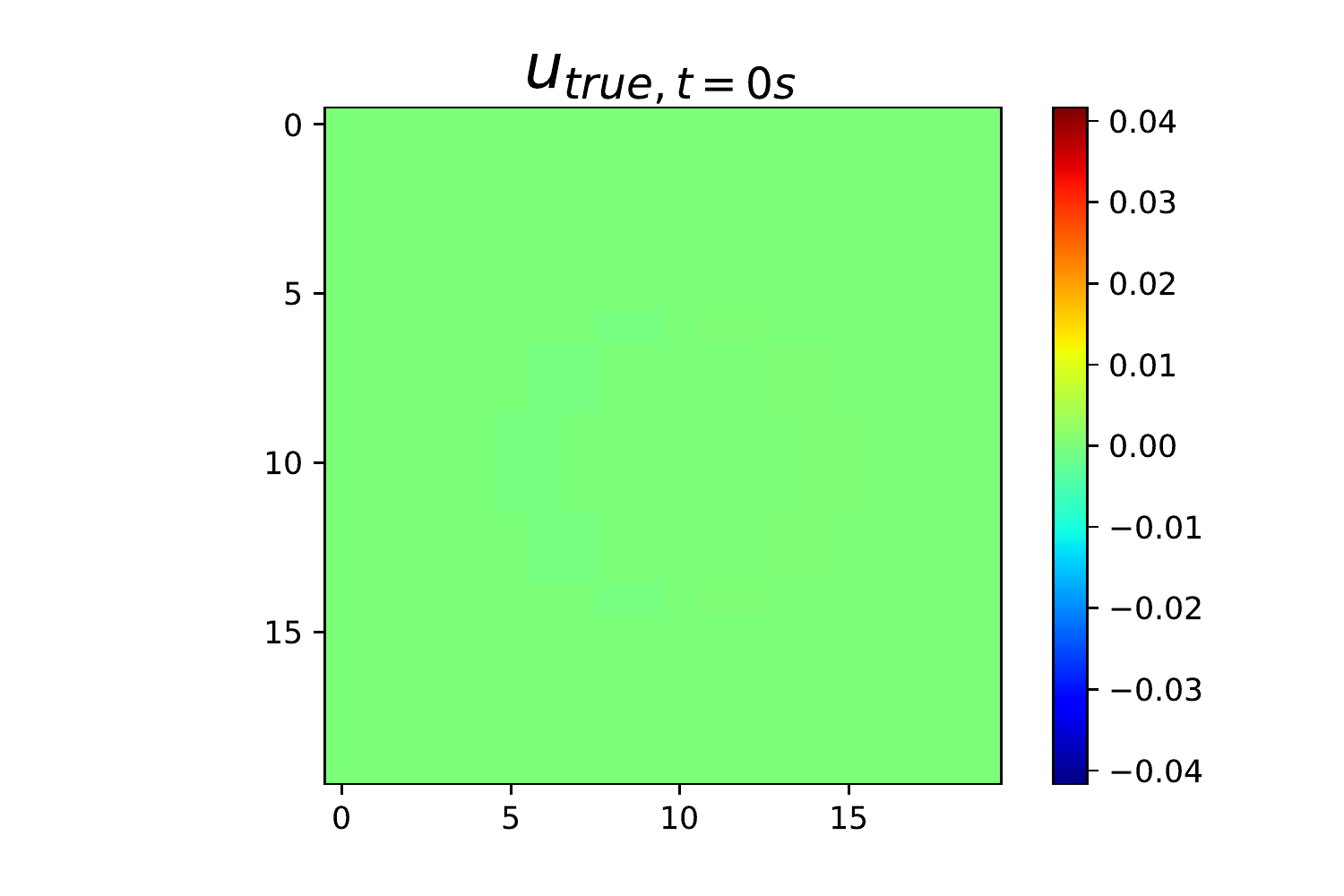}}
\subfloat[]{\includegraphics[width=2.2 in]{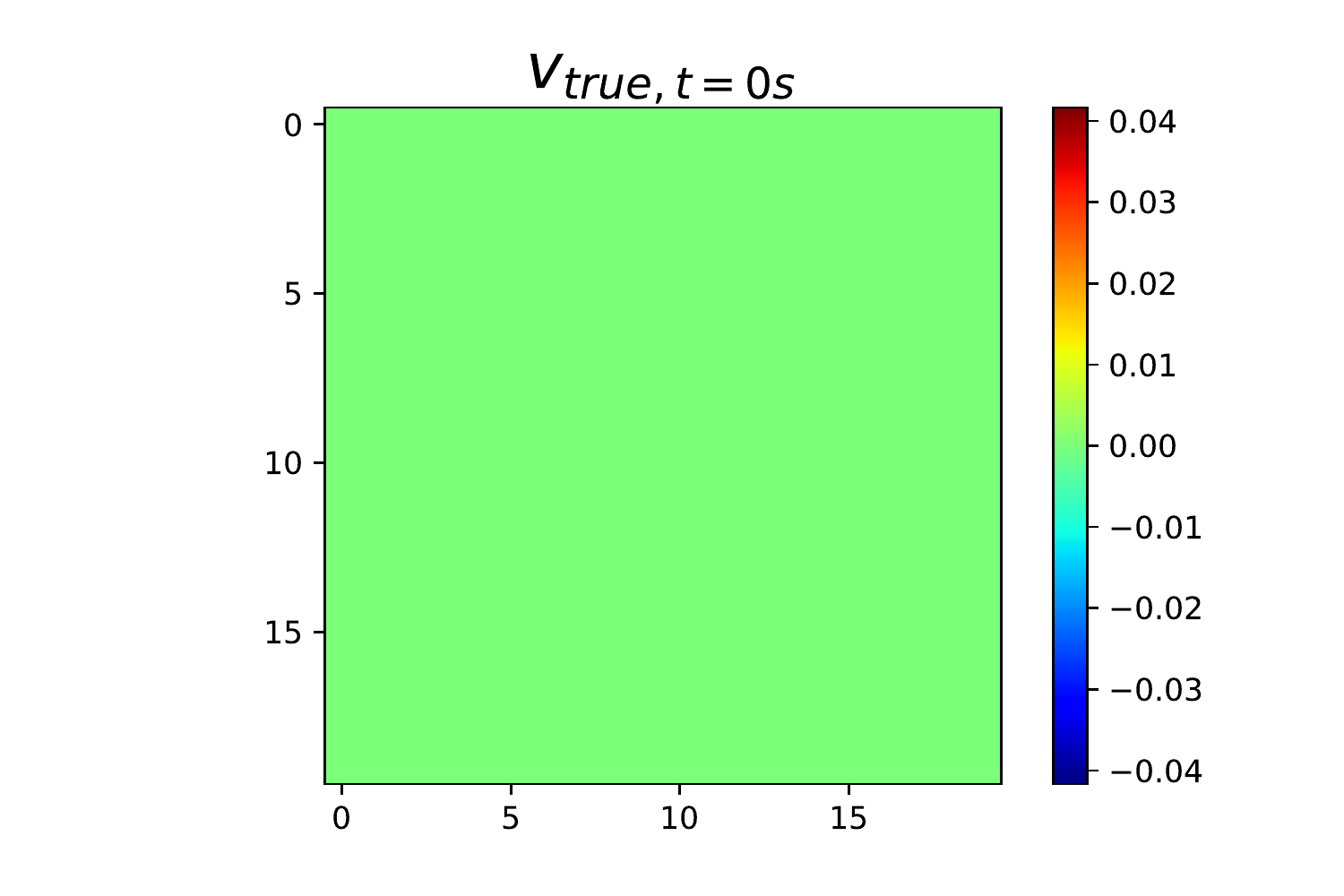}} \\
    \subfloat[]{\includegraphics[width=2.2 in]{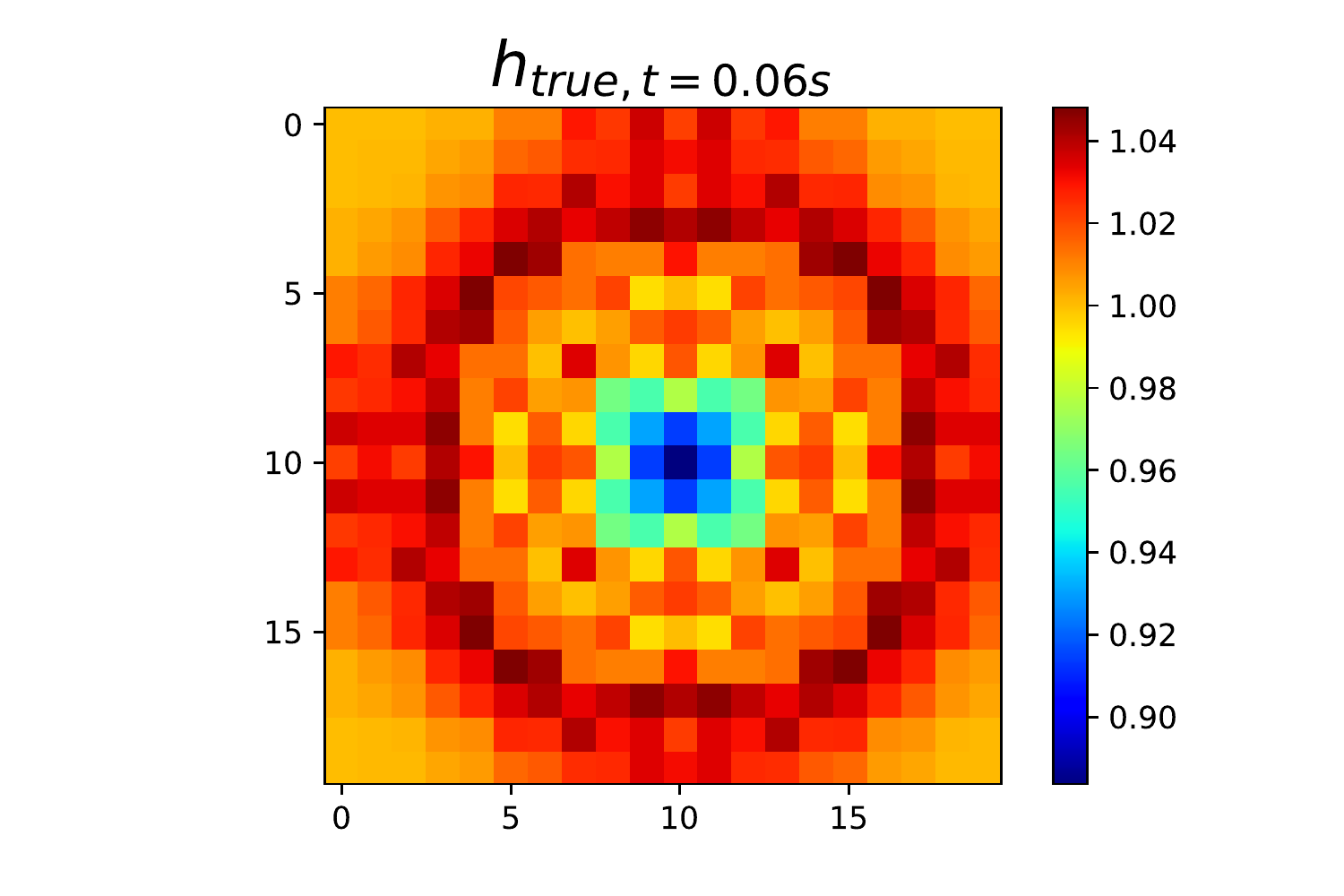}}
    \subfloat[]{\includegraphics[width=2.2 in]{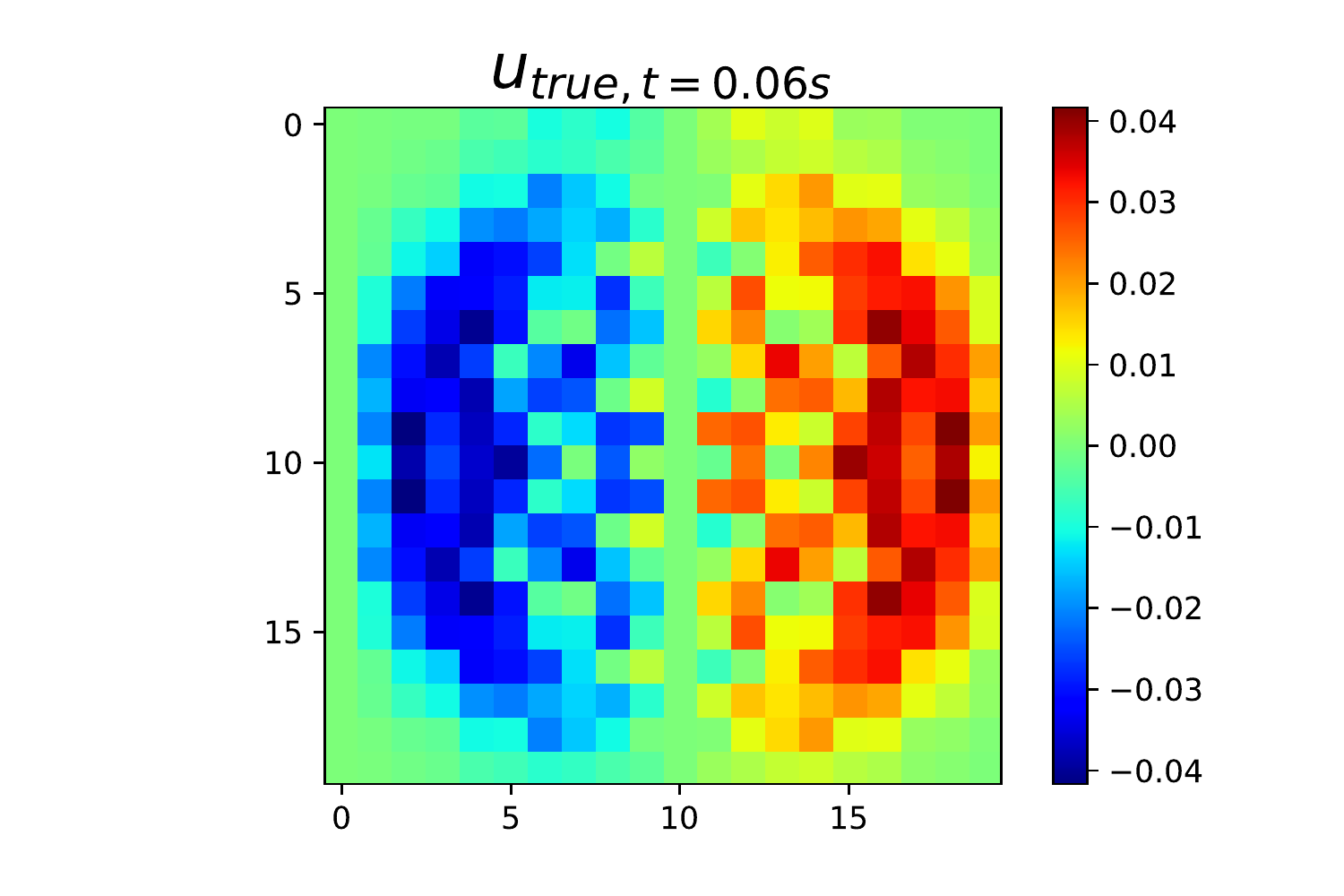}}
    \subfloat[]{\includegraphics[width=2.2 in]{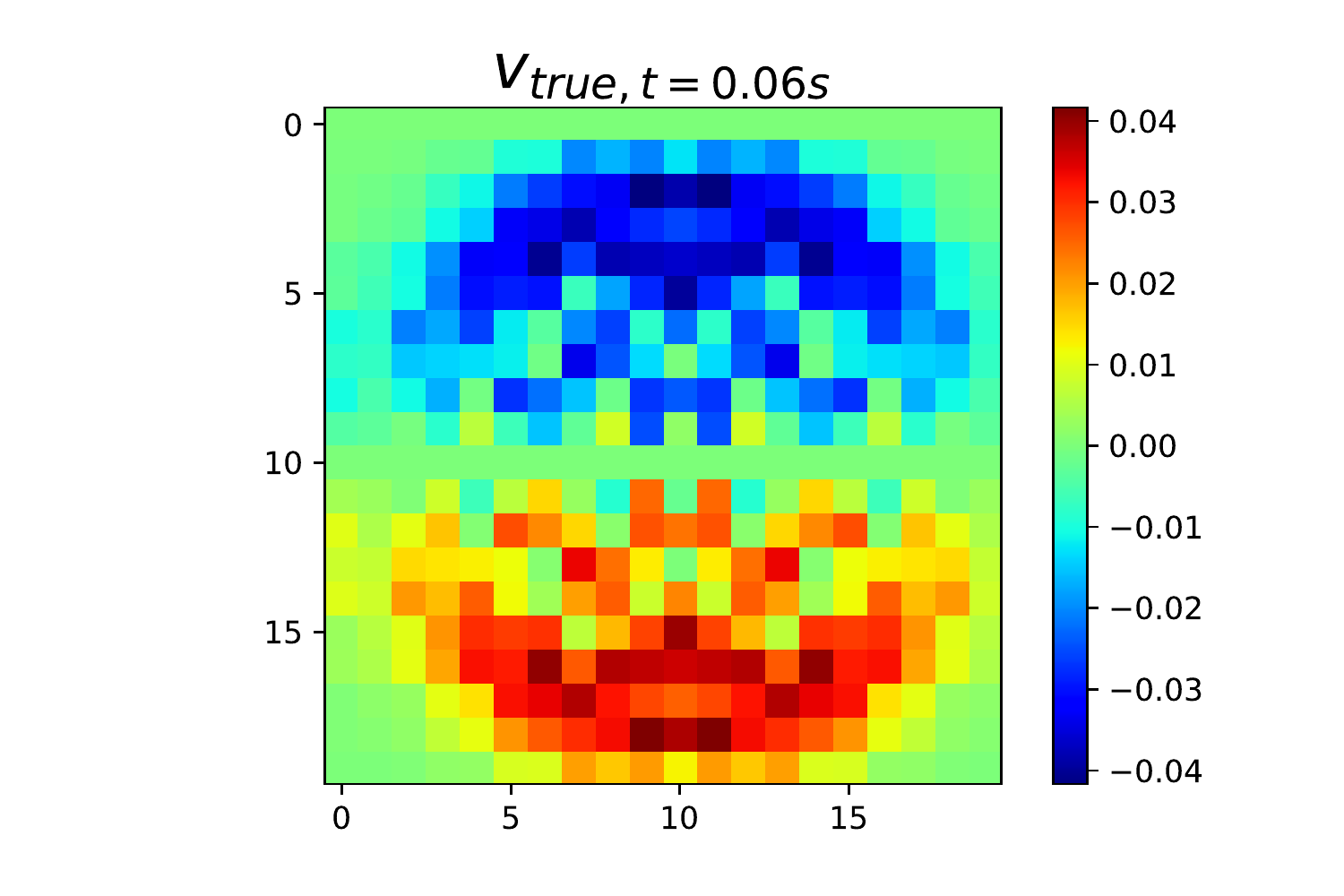}}\\
    \subfloat[]{\includegraphics[width=2.5 in]{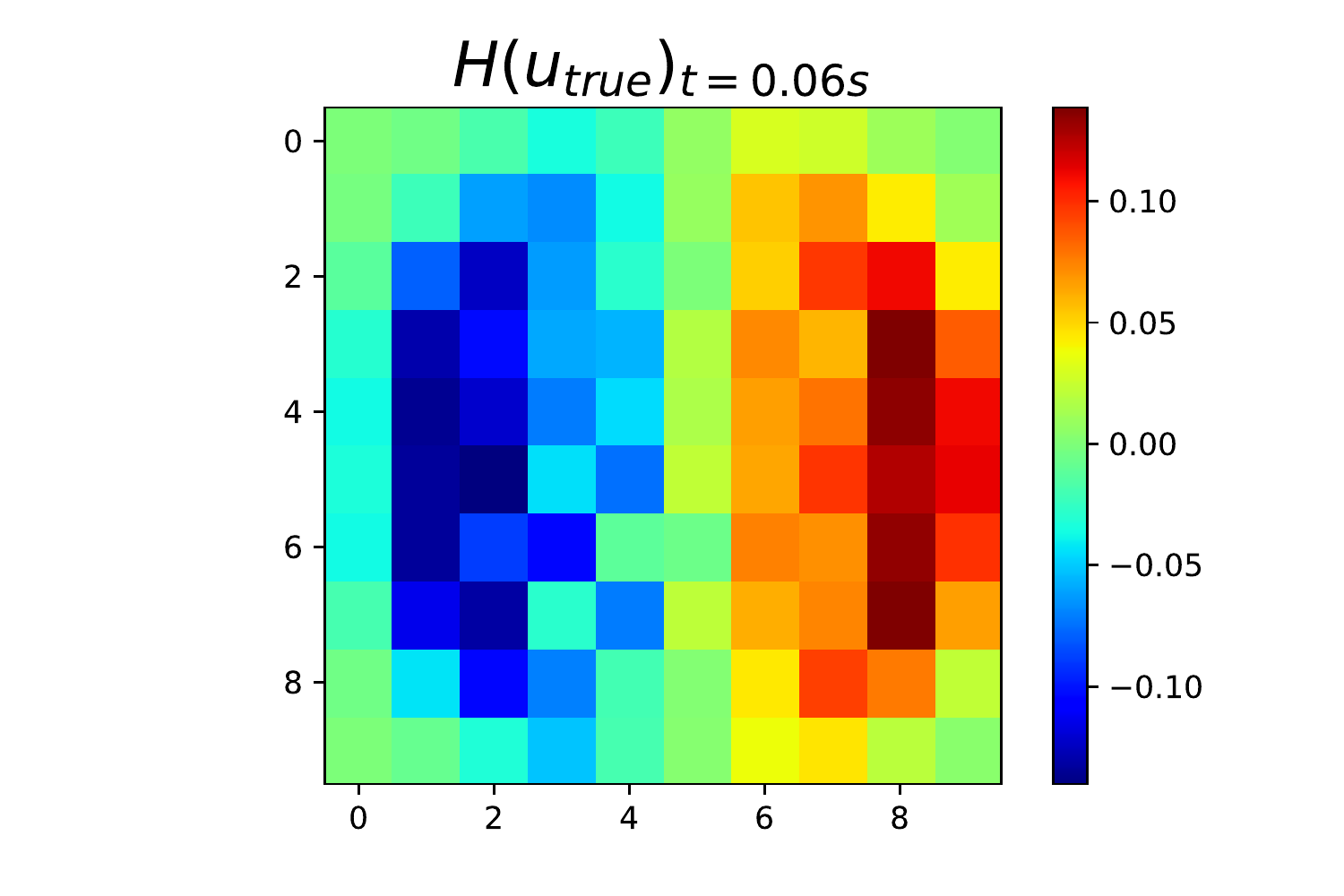}}
\subfloat[]{\includegraphics[width=2.5 in]{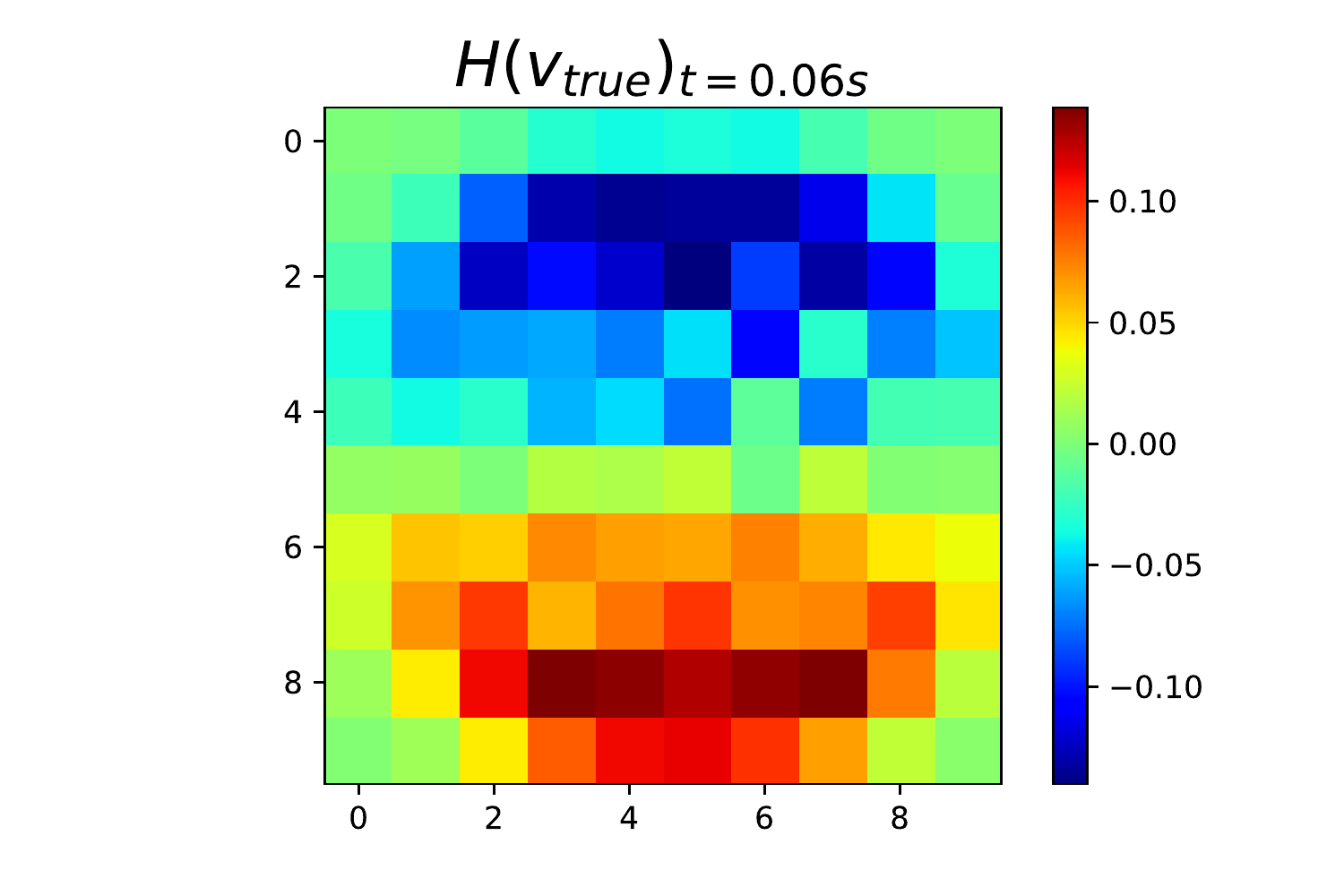}}
    \caption{Evolution of the shallow water model of $h,u,v$ (true states) at different time steps (a-f) and the error-free model equivalent $\bH(\bx_{true})$ for observations (g-h).}
  \label{fig:simulation}
\end{figure}

\vspace{4mm}
where $b=0.1$ is the viscous drag coefficient and the earth gravity constant $g$ is thus scaled to 1.
These equations are discretized in a $20 \times 20$ regular grid, solved by first-order finite difference method with a time discretization $\delta_t = 10^{-4}s$. This resolution is considered as the reference (i.e. the true state $\textbf{x}_\textrm{true}$) latter when performing DA algorithms. The state variables in this DA modeling are the combination of the velocity fields $\{u\}_{20 \times 20}$ and $\{v\}_{20 \times 20}$. The evolution of the reference ($\textbf{x}_{\textrm{true},t}$)  state is illustrated in Fig.~\ref{fig:simulation}. Spatially correlated prior error is then generated artificially for simulating the background state with a standard deviation $\sigma_{b,0} = 0.2$, i.e.
\begin{align}
    \textbf{x}_{b,t=0} \sim \mathcal{N} (\textbf{x}_{\textrm{true},t=0}, \textbf{B}_{t=0}) \quad \textrm{where} \quad \textbf{B}_{t=0} = {\sigma_{b,0}}^2 \textrm{corr}(\textbf{B}). \label{eq:Bt0}
\end{align}
The background error correlation matrix $\textrm{corr}(\textbf{B})$ is set to be isotropic (rotational invariant), following the second-order auto-aggressive (SOAR, also known as Balgovind) function,
\begin{align}
    \phi_\mathbf{B} (r) = \left ( 1+\frac{r}{L_\bB}\right ) \exp(-\frac{r}{L_\bB}),
\end{align}
where $r$ denotes the spatial distance and $L_\bB$ is the correlation scale length, fixed as $L_\bB = 4$ in this application. Being part of Matern kernels, the SOAR function is often used in DA for prior error correlation modeling \cite{cheng2019},\cite{Gong2020} thanks to its smoothness and good conditioning.
The simulation of $\textbf{x}_{b,t} = [{u}_{b,t}, {v}_{b,t}] $ via the same discretization of Eq.~\ref{eq: sw} (except the initial conditions) is used as background states at time $t$ in the DA modeling.
 For the knowledge of the exact\footnote{Here, by the term ``exact'', we refer to the covariance truly corresponding to the prior errors present in the background state, no matter the level of optimality of the chosen assimilation scheme.} background error covariance $\bB_{\textrm{E},t}$ at different time, $10^3$ background trajectories $\{\bx_{b,t}^{\gamma = 1...1000} \}$ are independently generated via Eq.~\ref{eq:Bt0}. This exact matrix, hidden for compression approaches, is seldomly used to evaluate the performance of DA algorithms with reduced observation. To simulate an industrial context, only 10 trajectories $\{\bx_{b,t}^{\gamma = 1...10} \}$ are used in the piecewise estimation of $\bH \bB \bH^T$ of a flow-independent window, making the ensemble size (10) much smaller than the problem dimension ($20 \times 20 = 400$).

\begin{figure}
  \centering
  \includegraphics[ width=5.5 in]{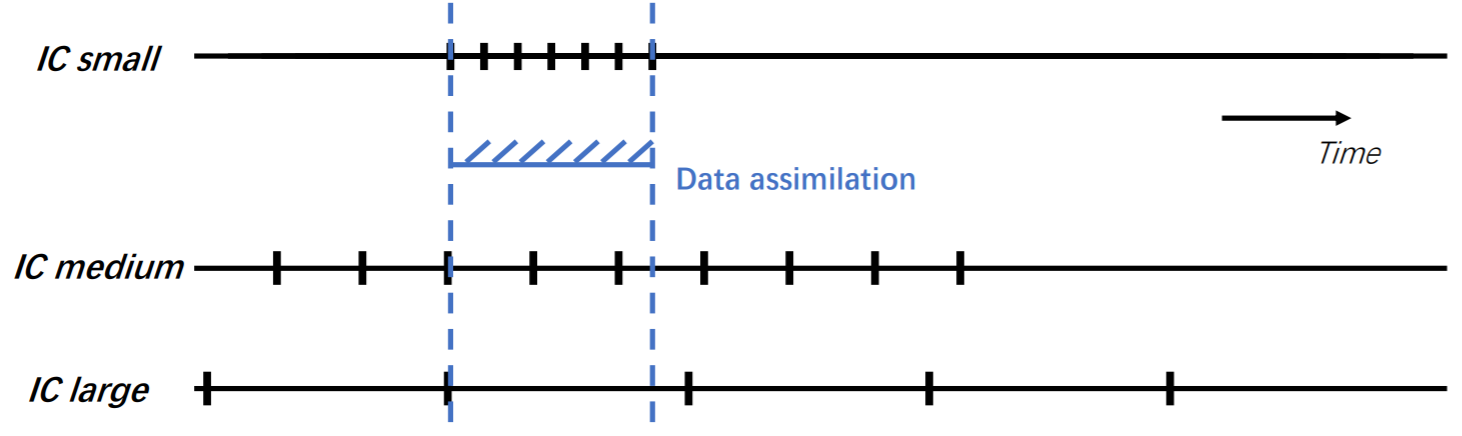}
    \caption{Simple sketch illustrating the three IC sampling strategies. The two vertical blue lines indicate where data assimilation experiments take place (as mentioned in Eq. \ref{eq:TrA}).}
  \label{fig:IC_samplings}
\end{figure}
The observations in these twin experiments are generated from the model equivalent based on the true states (i.e $\mathbf{H}(\bx_\textrm{true})$), separately for the fields $u$ and $v$, respectively denoted as $\textbf{y}_u$ and $\textbf{y}_v$.
For both fields, the observation $\by_{t} = [\by_{u,t}, \by_{v,t} ]$ at time $t$ is the sum of $u_t$ and $v_t$ in a $2 \times 2$ cells area with an observation error $\epsilon_{y_t}$,

\begin{align}
    \by_{u,i,j,t} = u_{\textrm{true},2i,2j,t} + u_{\textrm{true},2i+1,2j,t} + u_{\textrm{true},2i,2j+1,t} + u_{\textrm{true},2i+1,2j+1,t} + \epsilon_{y_{u,i,j,t}}
\end{align}

and identical for $\by_{v,i,j,t} $. Thus $\by$ represents also the evolution of the velocity field $u$ and $v$ with a "coarser" measure as shown in Fig.~\ref{fig:simulation} [g-h].

\begin{figure}
  \centering
  \subfloat[]{\includegraphics[width=2.9 in]{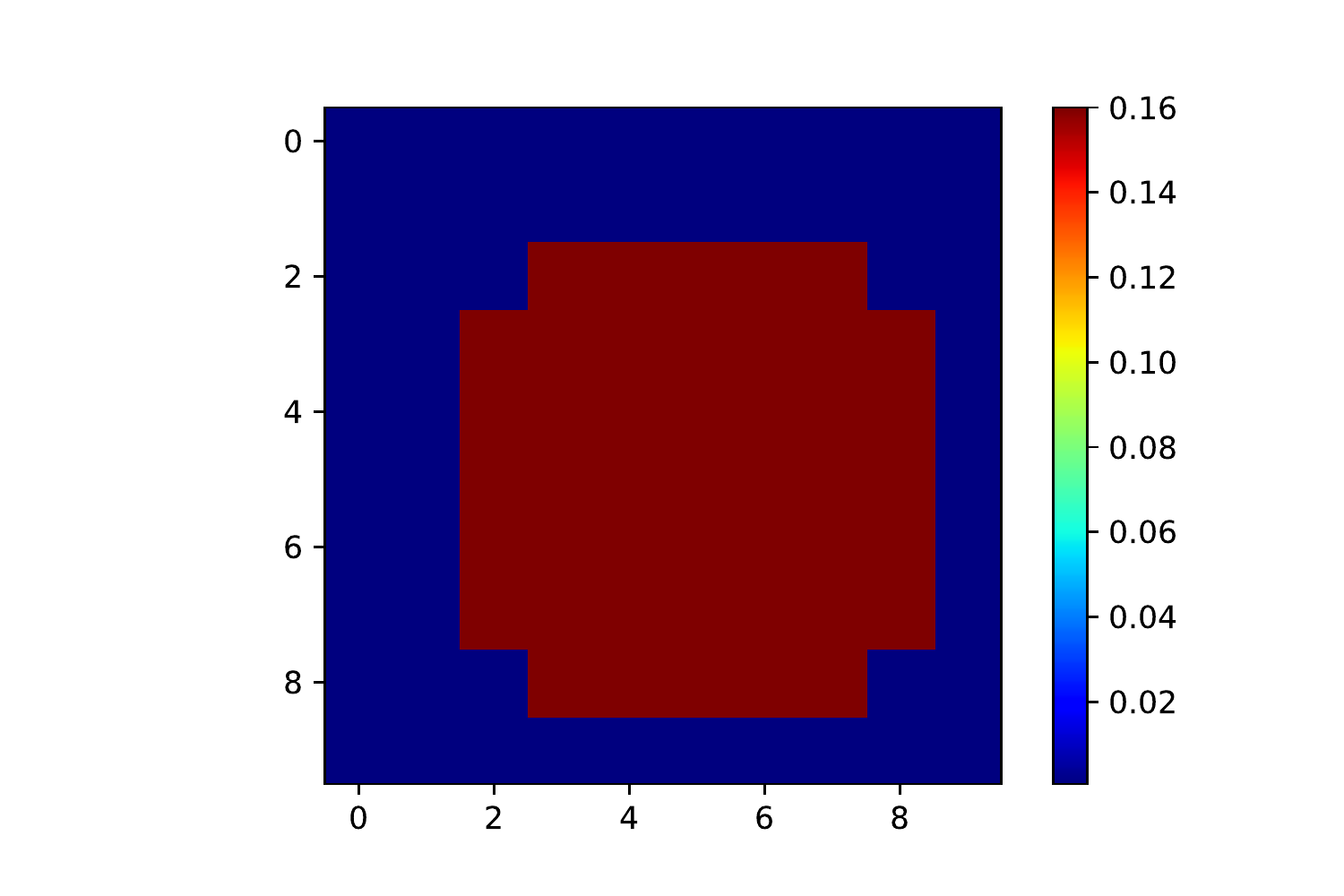}}
\subfloat[]{\includegraphics[width=2.9 in]{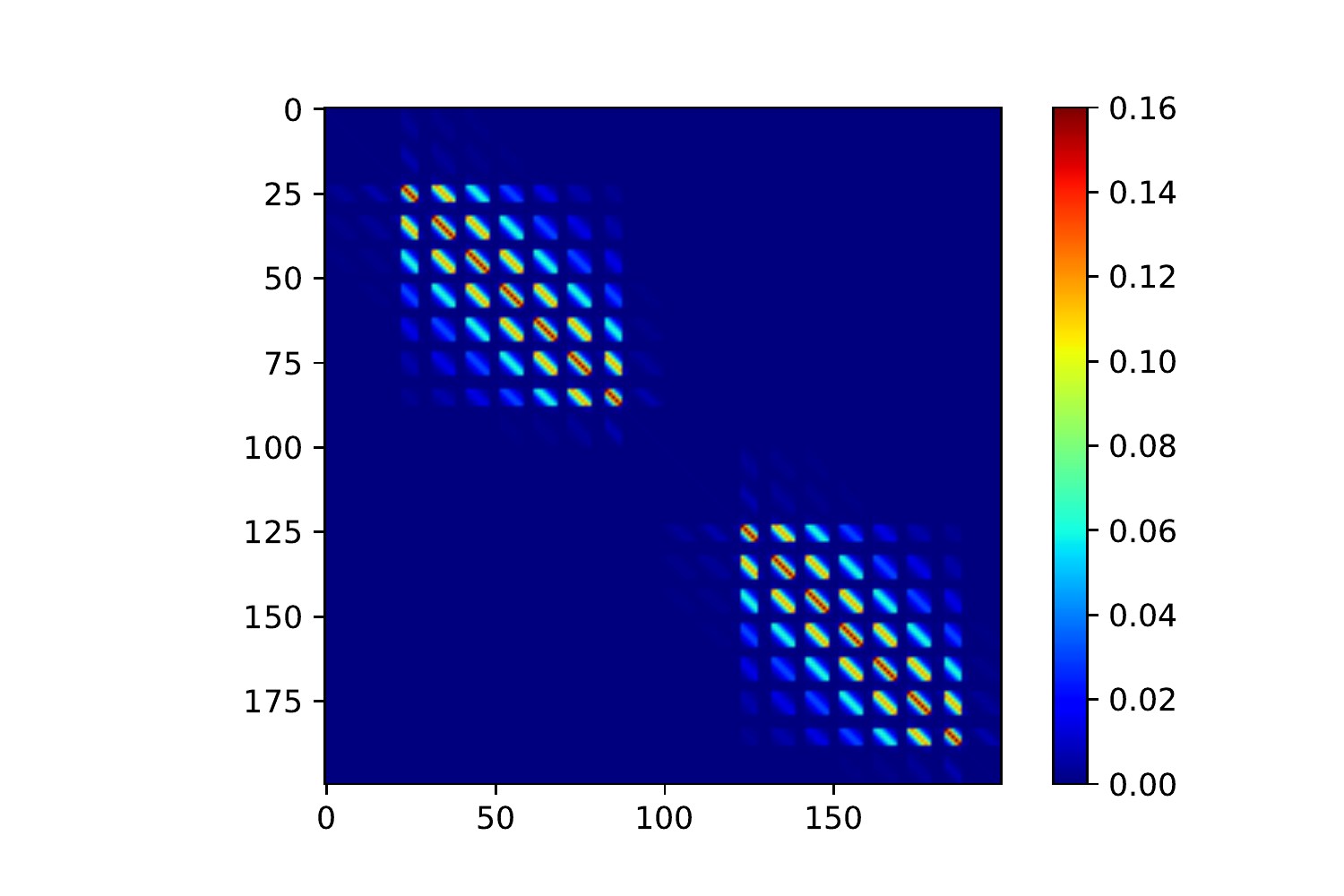}}
    \caption{The observation error variance of $\by_u$ and $\by_v$ in the shallow water model[a] and the Balgovind error covariances ($\bR$) after the observation vector (originally in a 2D grid) being converted to a 1D vector [b]}
  \label{fig:R_corr}
\end{figure}

In these experiments, we have set a non-homogeneous observation error covariance where the error deviation in the center (of radius $4$) of the field is 4 times higher, compared to boundary observations as show in Fig.~\ref{fig:R_corr}[a]. They are both of the same order of magnitude as $\sigma_{b,0}$, following also the SOAR function with a smaller scale length $L_\bR = 1$, compared to background error correlation. The full error covariance $\bR$ of observations $\mathbf{y}$ (after being converted to a 1D vector by concatenating rows of the original 2D grid model), supposed invariant against time, is illustrated in
 Fig.~\ref{fig:R_corr} [b]. The $\mathbf{R}$ matrix is supposed to be known in this application, thus only 10 observation trajectories $\{\by_{t}^{\gamma = 1...10} \}$ are generated to simulate an ensemble of small size while evaluating $\bH \bB \bH^T$ through Eq.~\ref{eq:hollings}. 
In this experiment, we make the choice to circumvent the difficulty
by setting a temporal correlated $\epsilon_{b,t}$, as the background noises are only added at the beginning of the simulation, and a temporal uncorrelated $\epsilon_{y,t}$. In fact, temporally correlated background errors are difficult to handle for the Desroziers method since it treats  the innovation quantities as independent samples for covariance estimation. These assumptions are realistic and widely adopted in DA problems since background simulations are often taken successively while observations are usually discrete. However, it is beneficial to have both time uncorrelated $\epsilon_{b,t}$ and $\epsilon_{y,t}$ for Desroziers-type estimation, as long as the error covariance could still be considered flow-independent \cite{Bathmann2018}.

\subsection{Numerical results for different compression strategies}
We then apply different strategies of observation compression and compare the performance of 3D-Var data assimilation using the reduced observation data. For each assimilation, only the current observation $\by_t$ is used to correct the background state $\bx_{b,t}$.
Thanks to the 1000 background trajectories $\{\bx_{b,t}^{\gamma = 1...1000} \}$ simulated, the exact $\bB_{\textrm{E},t}$ matrix can be empirically estimated at different time steps, allowing an accurate estimation of analysis error covariance $\bA_t$ via Eq.~\ref{eq:BLUE} since the $\bR$ matrix is supposed to be known. The matrix trace $\textrm{Tr}(\bA)$ then represents the sum of marginal analysis error, equivalent to the square of $L^2$ norm, i.e $\mathbb{E}\left( ||\bx_a-\bx_t||_2^2 \right)$, often used as an important indicator of DA schemes \cite{cheng2019}.

Another objective of this experiment is to inspect the impact on the assimilation error given by different sampling densities, which is critical in information-based compression, as stated in the introduction. 
  We display three sampling strategies  for $\bH \bB \bH^T$ estimation with different assumed flow-independent periods $[T_\textrm{s}, T_\textrm{f}]$,
\begin{itemize}
    \item \textit{IC small}: Dense sampling in a small period, $\Delta_t = 0.001s$ with $T_\textrm{s} = 0.16s$ and $T_\textrm{f} = 0.18s$
    \item \textit{IC large}: Sparse sampling in a long period, $\Delta_t = 0.1s$ with $T_\textrm{s} = 0s$ and $T_\textrm{f} = 2s$   
    \item \textit{IC medium}: Between \textit{IC small} and \textit{IC large}, with $\Delta_t = 0.01s$ $T_\textrm{s} = 0.1s$ and $T_\textrm{f} = 0.3s$,
\end{itemize}

as shown in Fig. \ref{fig:IC_samplings}, where $\Delta_t$ is the uniform time discretization between two snapshots. For all these three strategies, the $\bH \bB \bH^T$ is estimated via 20 time steps (i.e $T_\textrm{f}-T_\textrm{s} = 20 \Delta_t$), each with 10 background ($\{\bx_{b,t}^{\gamma = 1...10} \}$) and observation ($\{\by_{t}^{\gamma = 1...10} \}$) states/residuals. To gain a robust comparison, the posterior error variance $\mathcal{E}_\textrm{posterior}$ is averaged using $\textrm{Tr}(\bA_t)$ at four different time, included in all three assumed flow-independent windows,
\begin{align}
    \mathcal{E}_\textrm{posterior} = \frac{\sum \textrm{Tr}(\bA_t)}{4} \quad \textrm{for} \quad t \in \{ 0.16, 0.165, 0.170, 0.175\}. \label{eq:TrA}
\end{align}

\begin{figure}[t]
  \centering
  \subfloat[]{\includegraphics[width=0.3 \textwidth]{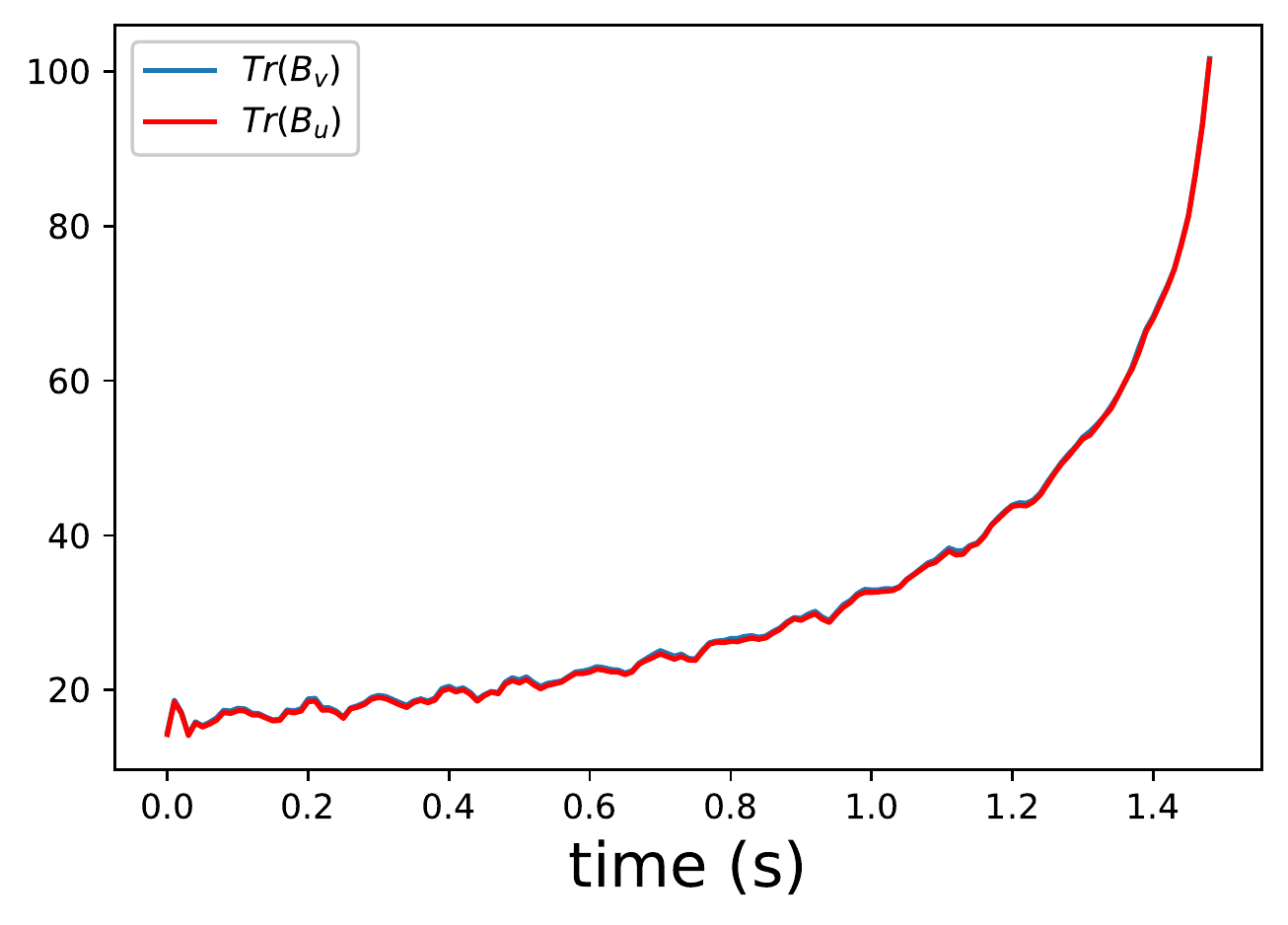}}
    \subfloat[]{\includegraphics[width=0.35 \textwidth]{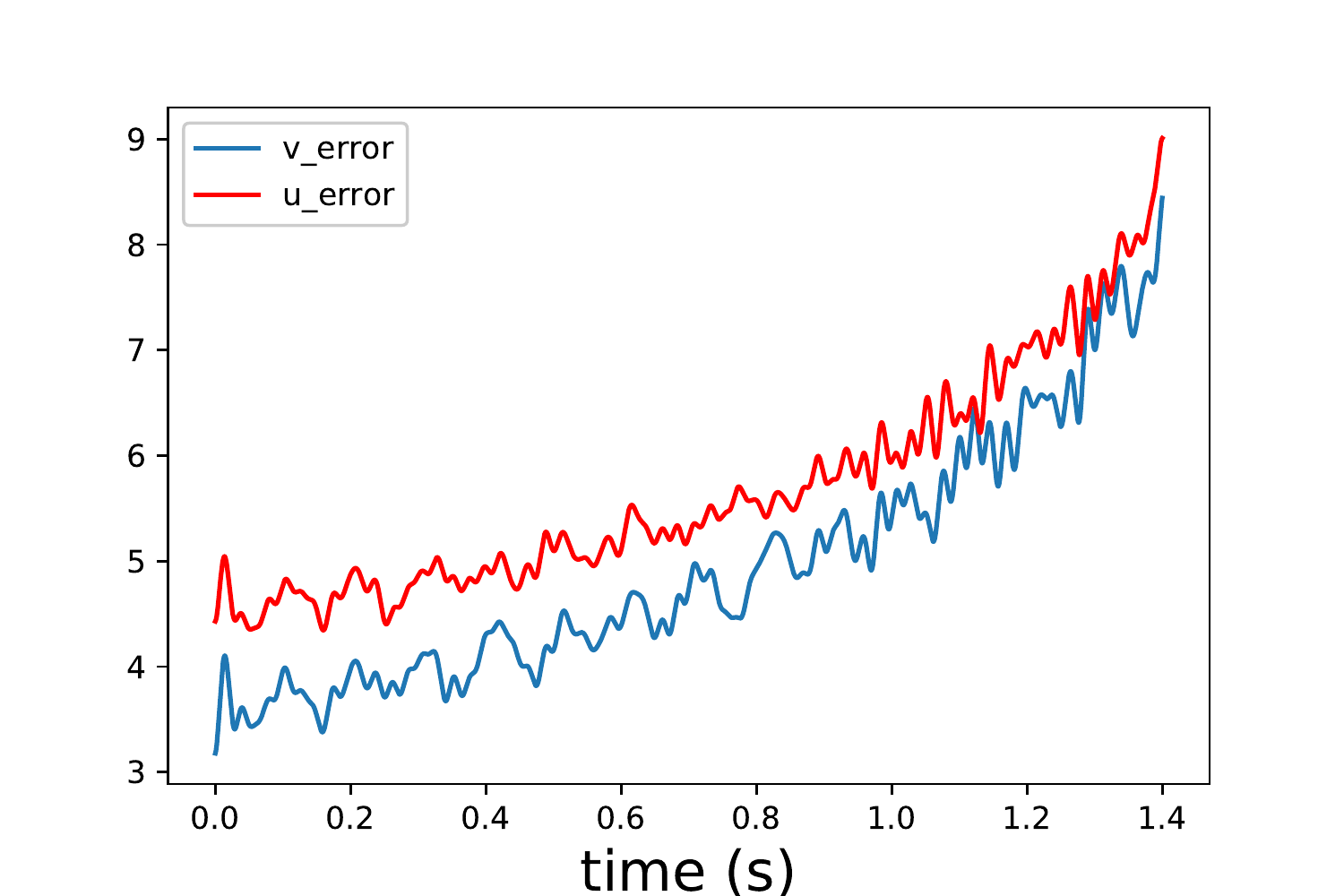}}
  \subfloat[]{\includegraphics[width=0.3 \textwidth]{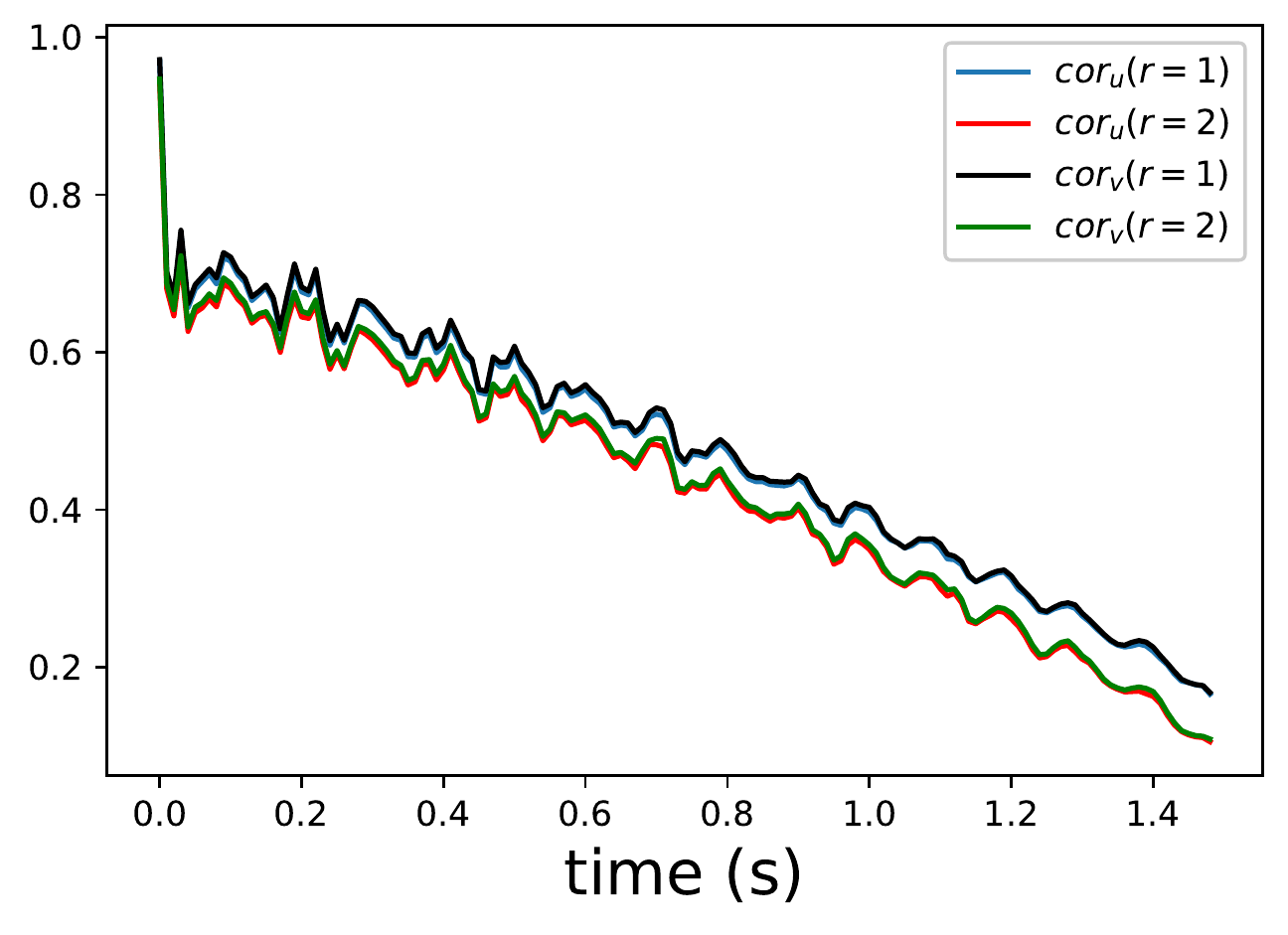}}\\
  \subfloat[]{\includegraphics[width=4.2 in]{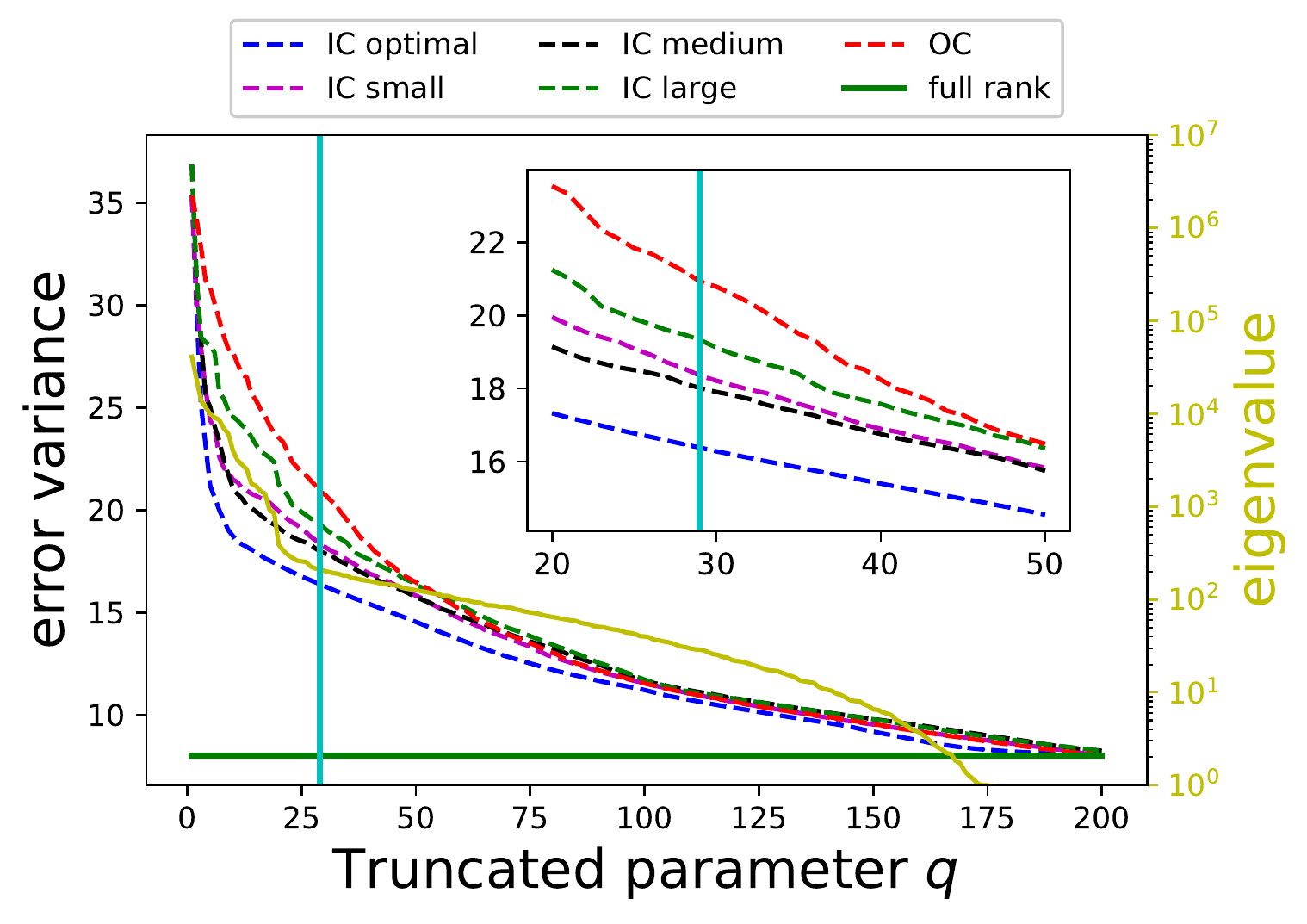}}
    \caption{ [a]: evolution of the exact background variance of $u$ ($\textrm{Tr}(\textbf{B}_u)$) and $v$ ($\textrm{Tr}(\textbf{B}_v)$) against time; [b]: evolution of $||\bx_{b,t}-\bx_{\textrm{true},t}||_2$ of a single background trajectory; [c]: evolution of average error correlation, of fixed distances ($r=1$ and $r=2$) in the 2D space; [d]: analysis error variance $\mathcal{E}_\textrm{posterior}$ (left y-axis) and eigenvalues of the estimated background error covariance in observation space ($\bH \bB \bH^T$) (right y-axis) of the medium sampling strategy as a function of the truncation parameter for $t\in [0.16s, 0.18s]$.} The vertical line represents the stopping criteria of $\sigma_{q} = \sqrt{\sigma_{1}}$.
  \label{fig:result}
\end{figure}

We illustrate in Fig.~\ref{fig:result}[d], the evolution of $\mathcal{E}_\textrm{posterior}$  against the truncation parameter $q$, varying  from 0 to 200. In fact, when $q=200$, all methods are equivalent since we work with the full observation data. From Fig.~\ref{fig:result}[d], we observe that all the information-based strategies with different sampling densities are always more optimal compared to the observation-based method for $q \in (0,200)$. We apply the stopping criteria as described in section \ref{sec:stop}, by calculating the eigenvalues of $\bH \bB \bH^T$ for the medium sampling strategy. We obtain the optimal truncation parameter $q_\textrm{optimal} = 29$. The distribution of these eigenvalues are shown by the right vertical axes in Fig.~\ref{fig:result}[d](numerical log scale is represented by the right vertical axis in matching color). With 29 modes, the assimilation correction is achieved from 53.3\% to 69.8\%, compared to the background model equivalent $\mathcal{H}(\bx_b)$ as shown in table \ref{table:2}, which is compatible to the results obtained in \cite{Fowler2019} when $L_\bB > L_\bR$.
Among the three sampling strategies, the one of "IC medium" owns the lowest output error variances, close to the optimal information-based compression where the $\bH \bB \bH^T$ is computed directly using $\bB_{\textrm{E},t}$. The latter, drawn with blue color in Fig.~\ref{fig:result}[d], stands for an optimal target for all information-based approaches since we suppose the exact background matrix is out of reach for data compression. As shown in this experiment, the choice of sampling strategy can significantly impact the compression optimality. If the samplings are too close, the residuals might not be uncorrelated, and if the samplings are too sparse, the flow independence of the $\bB$ matrix could be threatened. We remind that the stopping criteria for the truncation parameter $q$ varies for the different sampling strategies as shown in table~\ref{table:2}. However, in this experiment the values of the optimal truncation parameters obtained do not qualitatively change the results as shown in Fig.~\ref{fig:result}[d].

In Fig.~\ref{fig:result}[a,c], we display the evolution of the exact background error variances (i.e. $\textrm{Tr}(\bB_{\textrm{E},t})$) and  error correlation (for fixed distances, $r=1$ and $r=2$) against time. The estimation of background error correlation in the 2D space, also based on $\bB_{\textrm{E},t}$, is calibrated using the same method shown in \cite{cheng2019}.  We observe that the error variances increase continuously for both $u$ and $v$ while the spatial error correlation tends to shrink, both being significantly time-variant between $t=0s$ and $t=1.4s$. In order to illustrate the non-linear and turbulent nature of error propagation, we show in Fig.~\ref{fig:result}[b] the error evolution $||\bx_{b,t}-\bx_{\textrm{true},t}||_2$ of a single background trajectory. Obviously, these facts lead to problem of flow independent assumption for the \textit{IC large} approach (between $0s$ and $2s$), conducing a less optimal compression strategy as shown in Fig.~\ref{fig:result}[d]. From this twin experiment, we notice the advantage of information-based compression by selecting the most impacting observation components. The optimal sampling strategy may strongly depend on the characteristics (e.g chaosity, stability) of the dynamical system. \\

\begin{figure}[t]
  \centering
\includegraphics[width=2.5 in]{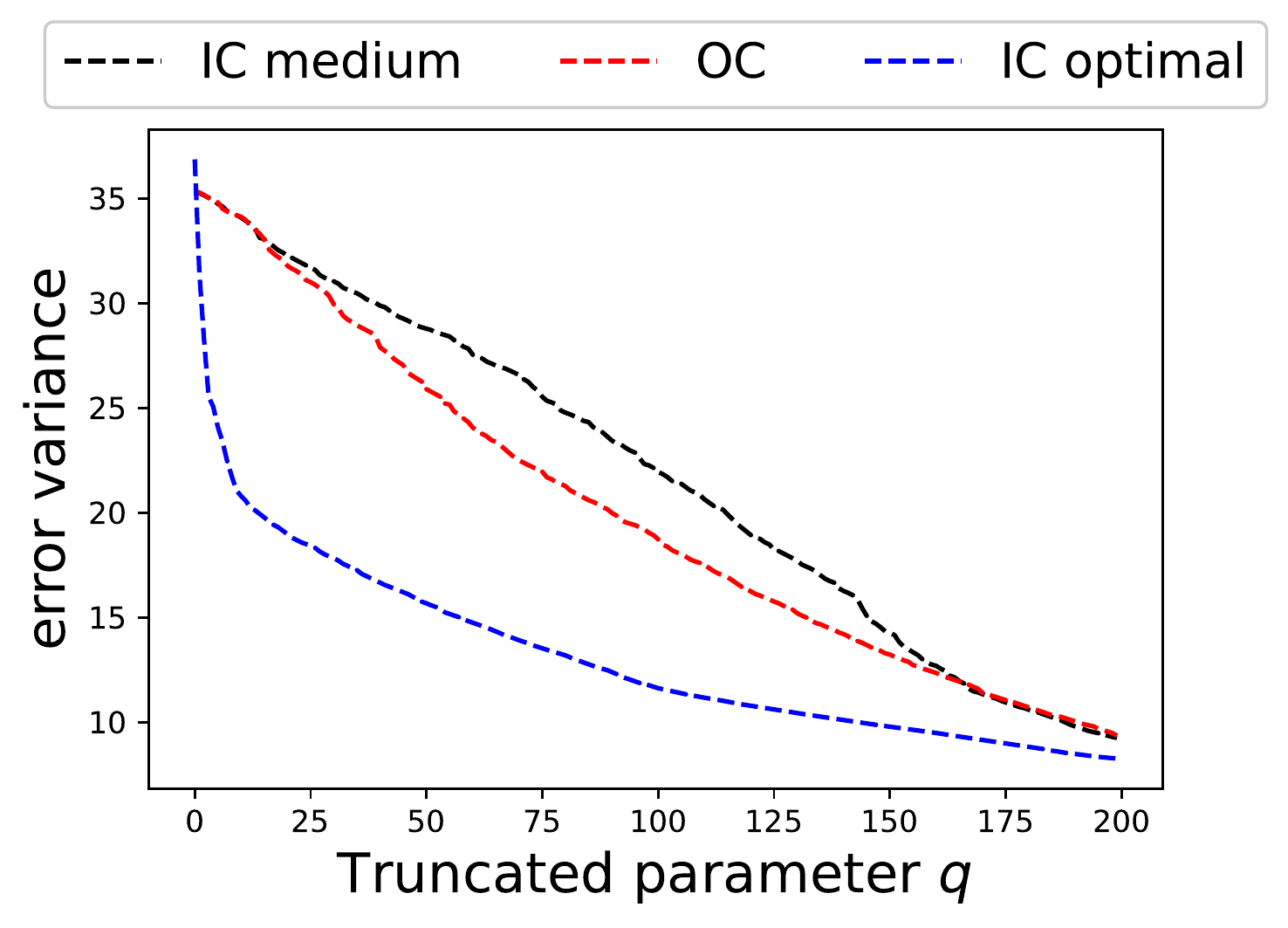}
\includegraphics[width=2.5 in]{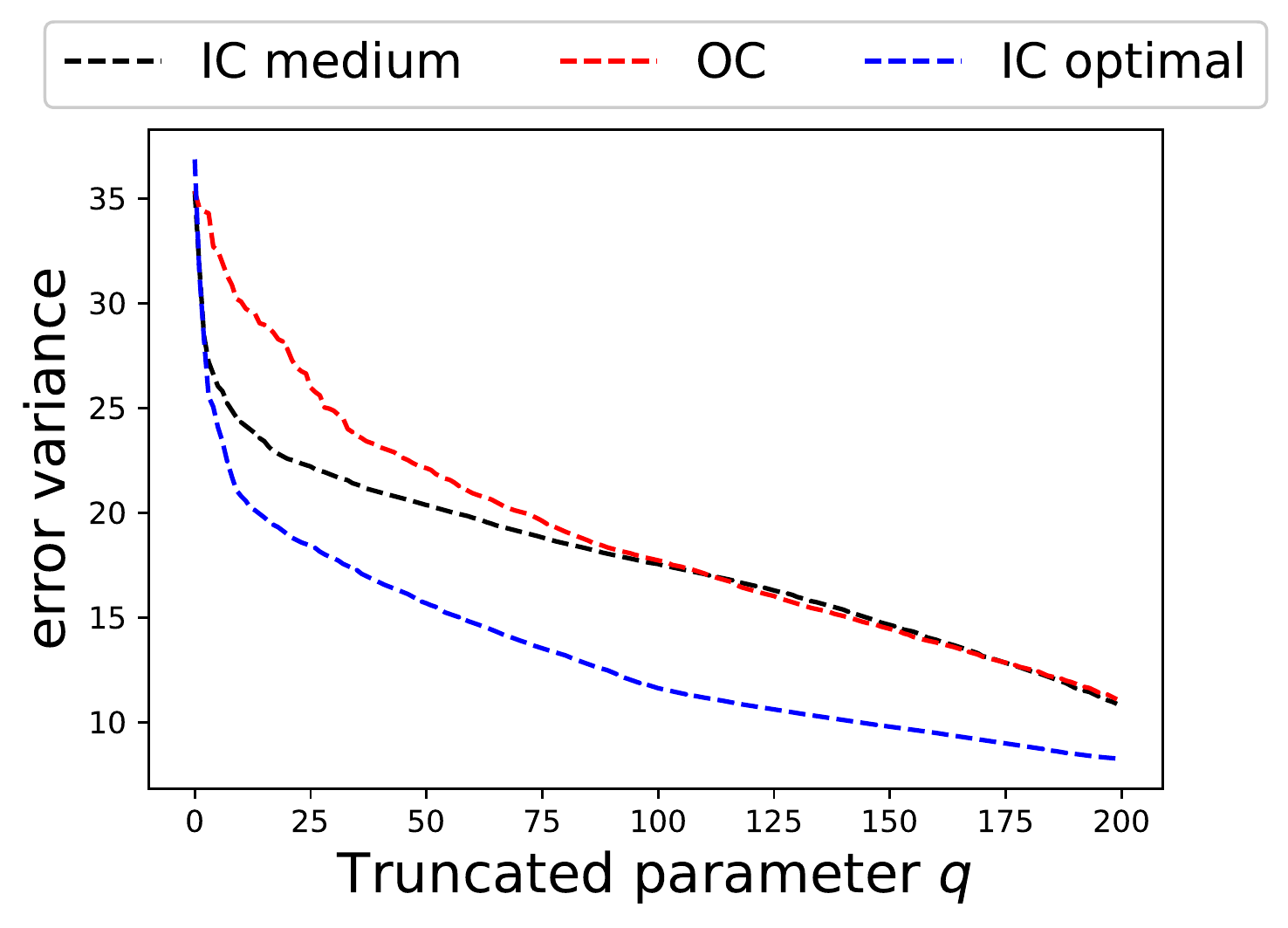}\\
(a) \hspace{5.5cm} (b)
    \caption{Evolution of error variance when the observation matrix is mis-specified.}
  \label{fig:mispec_R}
\end{figure}

Until now, we have shown that, in the idealised case where the observation matrix is known \textit{a priori} and the transformation operator is time-invariant, the information-based approach exhibits advantageous performance compared to the observation-based approach. However, as pointed out by \cite{Fowler2019}, IC approach can be sensitive to prior errors of covariance estimation. In order to investigate the impact of a potential misknowledge of matrix $\mathbb{R}$, we present here two cases where the difference between the assumed/estimated matrix $\bR_\textrm{A}$ and the exact matrix $\bR$ is voluntarily large. We explore two cases where the amplitude and the structure are misspecified, respectively:

\begin{itemize}
    \item (a): $\bR_\textrm{A}$ has the same correlation structure as $\bR$ with an homogeneous marginal error variance (i.e. $\bR_{\textrm{A},i,i} = 0.04$ which is different to $\bR$ (cf. Fig.~\ref{fig:R_corr}(a)).) 
    \item (b):  The correlation scale $L_{R_\textrm{A}}$ is set to be 5 while $L_{R} = 1$ as explained in section~\ref{sec:Experiments set up} with same marginal error variances.  
\end{itemize}
In both cases, the observation compression is implemented using $\bR_\textrm{A}$ while the observation matrix in the reduced space is set to be $\bR^{-1/2}_\textrm{A} \bR \hspace{1mm} \bR^{-1/2}_\textrm{A}$ instead of the identity matrix in Eq.~\ref{eq:OC_R} and Eq.~\ref{eq:IC_R}. The performance of these compression methods is illustrated in Fig.~\ref{fig:mispec_R}, respectively for case (a) and (b). The optimal IC solutions (same as the blue lines in Fig.~\ref{fig:result}(d)) are drawn in dashed blue lines for comparison purposes. Both OC and IC approaches exhibit less optimal performance compared to Fig.~\ref{fig:result}. Furthermore, as shown in Fig.~\ref{fig:mispec_R}(a), the IC method can be more sensitive to the mis-specification of the $\bR$ matrix amplitude, leading in this case to larger output error variances while IC behaves better than OC for misspecified R matrix correlation length.

\begin{table}[H]
\centering
\begin{tabular}{ ||c|p{2cm}|p{2cm}|p{2cm}|p{2cm}|p{2cm}||}
 \hline
 & OC & IC large & IC medium & IC small & IC optimal\\
 \hline
 \hline
$q_\textrm{optimal}$ & 22 & 48 & 29 & 25 & 78 \\
  \hline
Correction for $q=29$ & 53.3\% & 61.5\% & 65.7\% & 62.7\% & 69.8\% \\
  \hline
 \hline
\end{tabular}
\caption{The ratio of background minus analysis innovation ($||\mathcal{H}(\bx_b)-\mathcal{H}(\bx_a)||_2$) using compressed observation, relative to the one obtained with full observation}
\label{table:2}
\end{table}

\section{Application to an operating hydrological model}
\label{sec:hydro}
\subsection{DA modelling for flow reanalysis/prediction}
\begin{figure}[htbp]
  \centering
  \includegraphics[width=5.5in]{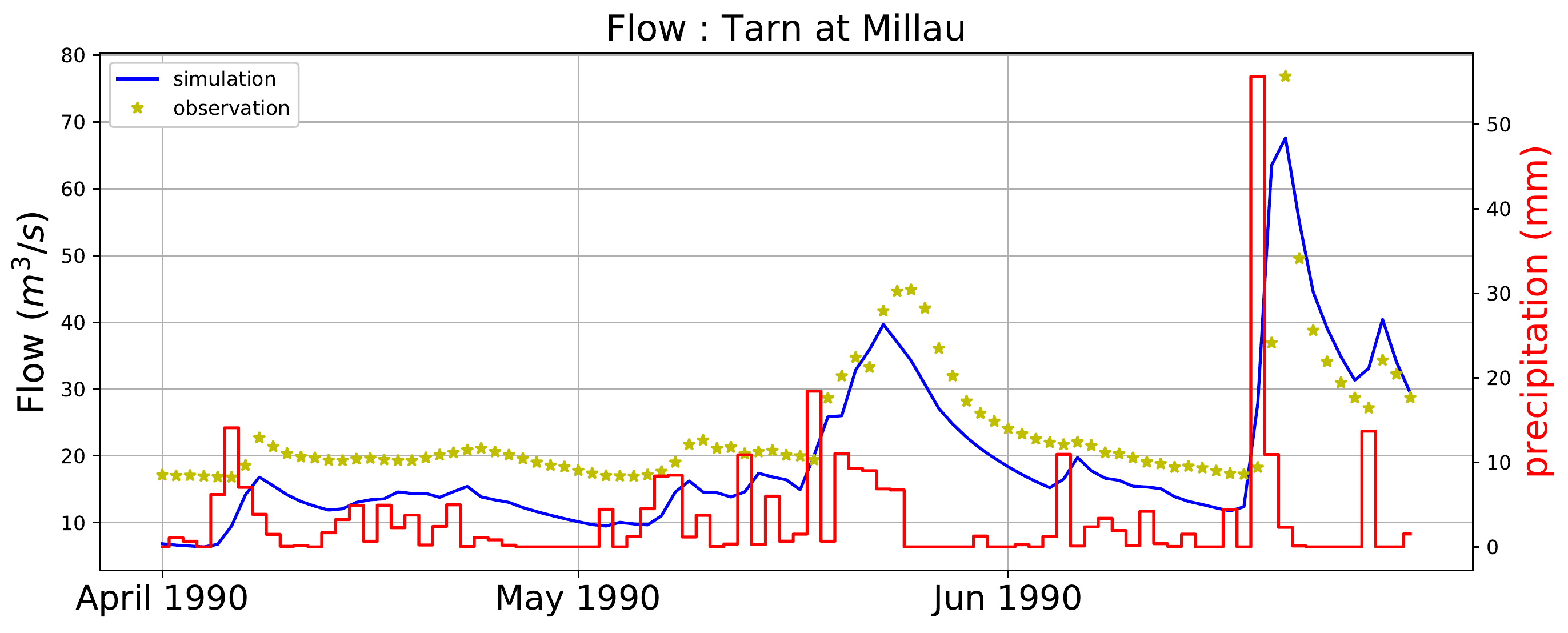}
     \caption{Example of simulation predicted by  MORDOR-TS using daily precipitation, and observed
     Tarn discharges at Millau for three months in 1990. Simultaneous observed
     precipitations are in red bars (with the scale on the right vertical
     axis).}
  \label{fig:mordor_ex}
\end{figure}

The compression strategies introduced in previous sections are applied to a hydrological application using a precipitation-flow simulator MORDOR-TS developed by Électricté de France (EDF, the French electric utility company). This software is widely applied in operating hydraulic/horological problems, e.g. \cite{garcon1996}, \cite{Rouhier2017}, \cite{chengthesis}. Based on information on spatially distributed physical parameters, such as precipitation or temperature, it provides a simulation of river flow relying on conceptual watersheds modeling. For more details about MORDOR-TS, interested readers are referred to \cite{Rouhier2017} and \cite{cheng2020b}. MORDOR-TS is used as a non-linear state-observation transformation operator in data assimilation. We concentrate on a study area in the south of France, around the Tarn river where 9 streamflow gauges positioned at different mesh outlets are available. The Tarn river, being known for its extreme variability of water-level values and high sensitivity to
precipitations \cite{cheng2020b}, is an ideal benchmark for
comparing different DA strategies. Located downstream, the Tarn river outlet at Millau (hereby denoted as TM) is of particular interest in the hydrological study. As an example, we show in Fig.~\ref{fig:mordor_ex} the simulated and daily observed
Tarn river discharges at Millau, for 3 months in 1990 with the
averaged precipitation over 28 spatially distributed regions (see \cite{cheng2020b}). Significant impacts of precipitation on the river flow of TM is observed with a delay of 2 to 5 days. The objective of this DA modeling is to improve the river flow  prediction and reanalysis (history matching) by performing corrections on the daily precipitation in the 28 regions. Other physical quantities (e.g temperature) are considered as invariant parameters in this study. The variational assimilation is performed using the ADAO \cite{adaotrue} package of SALOME platform, also developed by EDF. 

As mentioned in \cite{cheng2020b}, performing DA correction on all precipitation inputs (i.e 28 regions) can probably introduce an over-parameterization and thus induces an overfitting, with a high risk to deteriorate flow forecasts. Therefore, we make the choice to proceed with uniform additional increments $\xi^p_{t}$ for all 28 regions, depending only on time $t$. Incremental variables $\xi^{r,j} (j=1..8)$ on the eight parameters which determine the initial (at $t=0$) reservoir level is also added in the state space to adjust the river flow at the beginning of each assimilation window.
These windows are fixed of 30 days, leading to an observation vector of dimension 270 with 9 gauges. Temporal correlation is considered  for both background and observation errors. The DA modelling is summarized in Table~\ref{table:1} and a more detailed description can be found in \cite{cheng2020b} and \cite{chengthesis}. The main objective of this application stands for improving short-range flow forecasting by correcting historical precipitation. Since the impact of the precipitation on the river flow is only significant within 3 to 4 days (see \cite{cheng2020b} for details), we fix the prediction window to 3 days in this study.
\begin{table}[H]
\centering
\begin{tabular}{ ||c|c|c|c|c|c||}
 \hline
DA modelling & state: $\bx$ & $dim(\bx)$ & Observations: $\by$ & $dim(\by)$ & invariant parameters\\
 \hline
 \hline
Incremental 3DVar & \makecell{ $\xi^p_{t}$ \\  $\xi^{r,j}$} & $38$ & river flow {$Q_{q,t}$} & $270$ & { temperature, etc} \\
  \hline
 \hline
\end{tabular}
\caption{Details of DA modelling where $t= 0..29$ is the time (days) relative to the beginning of the assimilation window; $q= 1..9$ represents the 9 gauges where $\xi^p_{t},\xi^{r,j}$ represent respectively the increments of daily precipitation and initial reservoir level. }
\label{table:1}
\end{table}

\subsection{Observation compression}
Despite that MORDOR-TS is computationally efficient (it may take only a few CPU seconds to simulate a spatially distributed
flow simulation of several years), the application of variational assimilation algorithms could be expensive, due to the non-linearity of the transformation operator. As shown in Table. \ref{table:1}, in this DA modeling, the dimension of the observation vector is much larger compared to the state dimension, promoting the utilization of observation data compression.  
We then implement DA algorithms in the hydrological model with compressed data using either OC or IC approaches. To make the compression strategies more general, in both cases the principal components are constructed using the daily observed flow data from 1990 to 2000 in the 9 gauges. The objective of this study is to make an efficient use of the observation vector $\by$ with an optimal number of modes selected, which we expect to be much smaller than the full observation dimension ($\textrm{dim}(\by)=270$). 

A major hurdle of this application is that the \textit{a priori} knowledge of both $\bB$ and $\bR$ is very limited. As a remedy, we start as described in \cite{cheng2020b}, by considering the background covariance matrix $\bB$ of Balgovind-type since we wish to model the existence of temporal correlation in the precipitation data. Moreover, the initial $\bR$ matrix is set to be diagonal. The DI01 algorithm \cite{Desroziers01} is then applied several times to to come up with a reasonable approximation of the ratio between $\textrm{Tr}(\bB)$ and $\textrm{Tr}(\bR)$ at the first stage. In a second stage, we then perform the estimation of $\bH \bB \bH^T$ and $\bR$, relying on Desroziers formulation (respectively Eq.~\ref{eq:d05R} and \ref{eq:d05}) using 3400 assimilation windows of 30 days from 1990 to 2000. By then, post-processing is required to ensure the symmetric positive definiteness (SPD) of the $\bR$ matrix. More precisely,
\begin{align}
    \mathbf{R} \longleftarrow \frac{1}{2}  (1-\mu)(\mathbf{R} + \mathbf{R}^T) + \mu \textbf{C}, \label{eq:hybridSPD}
\end{align}
where $\mu = 0.1$ and $\textbf{C} = \textrm{Tr}(\mathbf{R}) \times
\textbf{I}$. The Desroziers method is iterated twice, using the same data set, to ensure the stability of the estimated matrices. The algorithm outputs produced after the first  and the second iterations are very similar as shown in \cite{cheng2020b}.
We emphasize that the estimated $\bR$ matrix is not only used for the observation compression but also in the DA algorithm in the full observation space. The $ \bH \bB \bH^T$ matrix is obtained through Eq.~\ref{eq:hollings}, once the $\bR$ matrix is specified. As a remark, even if the system considered here is not very large, the computational burden associated with the data assimilation of this nonlinear system (for which prior information is degraded) remains important because of a multi-stage tuning approach which combined several offline and online covariance tuning algorithms can be implemented to improve the reanalysis and the forecasting accuracy of this hydrological application. However, these methods are computationally expensive, especially when iterations are needed (e.g. \cite{cheng2019}). With advanced data compression methods, the computational burden can be released, allowing more precise covariance tuning to improve the DA performance. 

\subsection{DA with compressed data}
\subsubsection{Averaged performance}
\begin{figure}[H]
  \centering
    \includegraphics[width = 3. in]{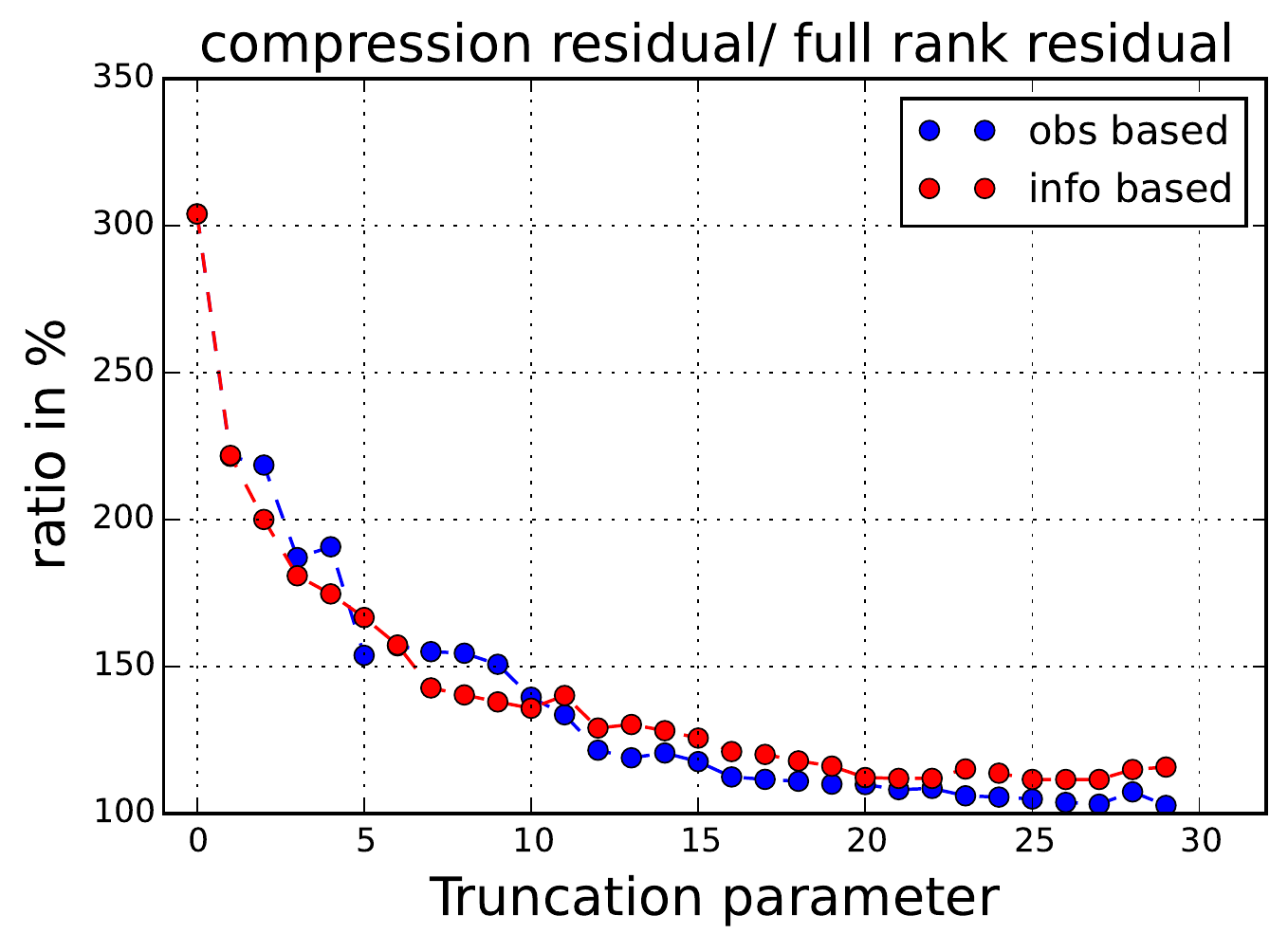}
    \includegraphics[width = 3. in]{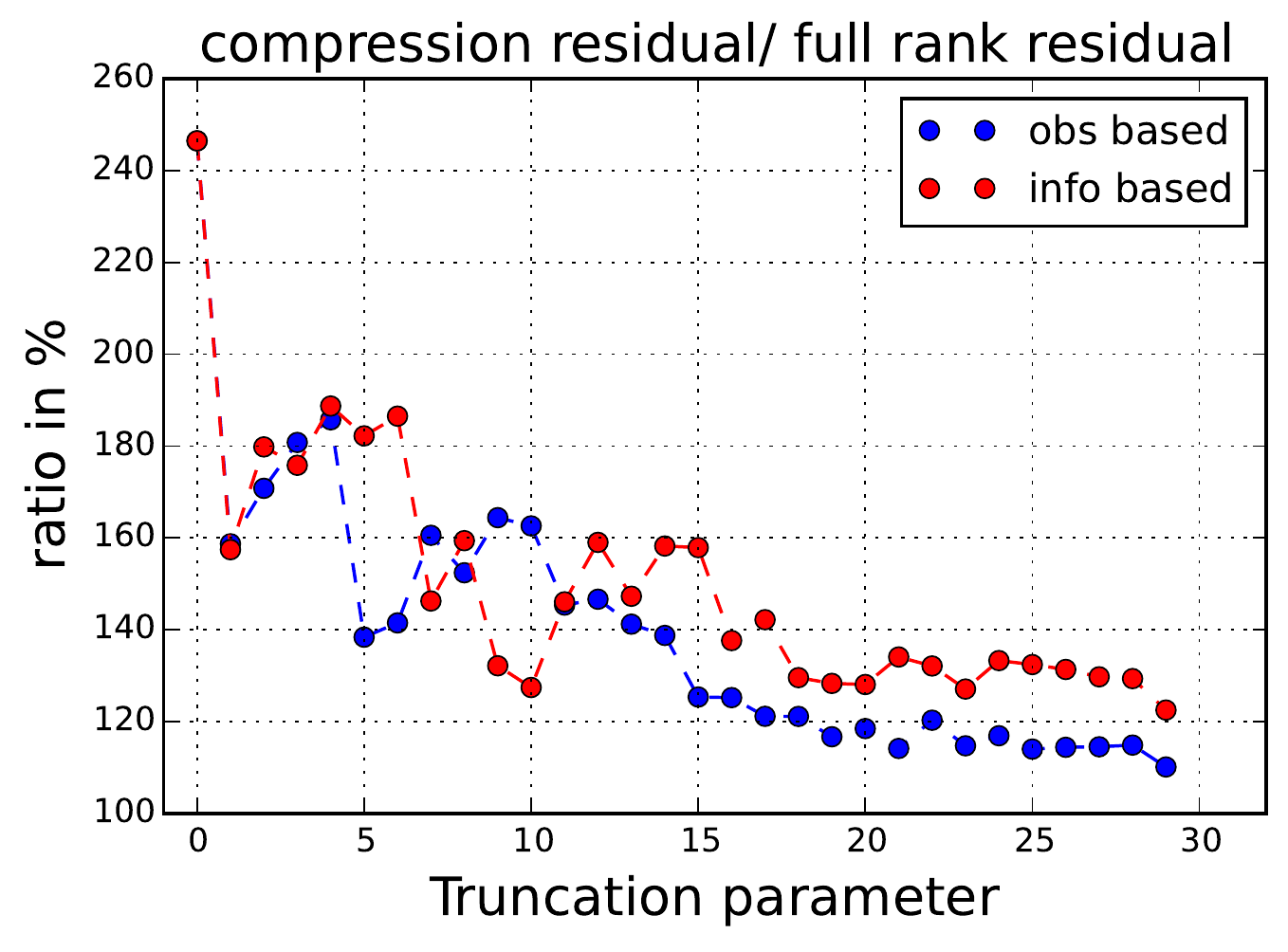}\\
    (a) \hspace{7.6cm} (b)
    \caption{ Evolution of $\mathcal{r}$, averaged using 12 assimilation windows, against the number of truncation parameter $q$ for reanalysis[a] and prediction[b] at TM.  }
  \label{fig:1193-1-1 flow}
\end{figure}

Extracting the principal components $\tilde{\bL}$ and $\hat{\bL}$, respectively based on estimated $\bH \bB \bH^T$ and $\bR$, we then apply the compression methodology described in sect.\ref{sec: metho}. The objective is to compare the assimilation output $\bx_{a,\textrm{compression}}$ and $\bx_{a,\textrm{full}}$,  obtained using either the compressed observation $\hat{\by}_q$, $\tilde{\by}_q$ or the full observation vector $\by$. More precisely, we are interested in the observation minus analysis (O-A) innovation quantity for both flow reanalysis and forecast. Varying the truncated parameter $q$, DA processes are performed respectively with  $(\bx_b, \tilde{\by}_q, \bB, \tilde{\bR}_q, \tilde{\mathcal{H}}_q)$ and  $(\bx_b, \hat{\by}_q, \bB, \hat{\bR}_q, \hat{\mathcal{H}}_q)$ for 12 assimilation windows in 1993, each of 30 days starting at the first day of every month. We draw the averaged compressed/full O-A innovation ratio $\mathcal{r}$, defined as
\begin{align}
    \mathcal{r} = \frac{||\by-\mathcal{H}(\bx_{a,\textrm{compression}})||_2}{||\by-\mathcal{H}(\bx_{a,\textrm{full}})||_2},
\end{align}
in Fig.~\ref{fig:1193-1-1 flow} for both reanalysis[a] and prediction[b] at TM. More particularly, $\mathcal{r} = 100\%$ means the reanalysis/prediction accuracy of the current solution is equivalent to the one obtained with the full observation vector.

We observe from Fig.~\ref{fig:1193-1-1 flow} that the performance of these two approaches is similar to the reanalysis while the observation-based method is slightly more optimal on average for flow forecasting. The evolution of the reconstruction error (Fig.~\ref{fig:1193-1-1 flow}[a]) is much smoother, compared to the prediction error ((Fig.~\ref{fig:1193-1-1 flow}[b])), both against the truncation parameter. In fact, the reconstruction error is estimated using assimilation windows of 30 days while prediction windows are solely of 3 days. Therefore, the estimation of the prediction ratio has significantly more sampling noise. Furthermore, since both $\bB$ and $\bR$ are not well specified \textit{a priori}, extra noise can be introduced while estimating the information entropy. For both methods, the assimilation results obtained using 15 to 20 modes (around $5\%$ to $7.5\%$ of total observation dimension) are close to the full rank solution in terms of both reanalysis and prediction. Without deteriorating the assimilation result, these compression strategies make the DA algorithm certainly more efficient, allowing more optimization iterations if needed.

\subsubsection{Performance in each DA window}
\begin{figure}[htbp]
  \centering
    \includegraphics[width = 2.2 in]{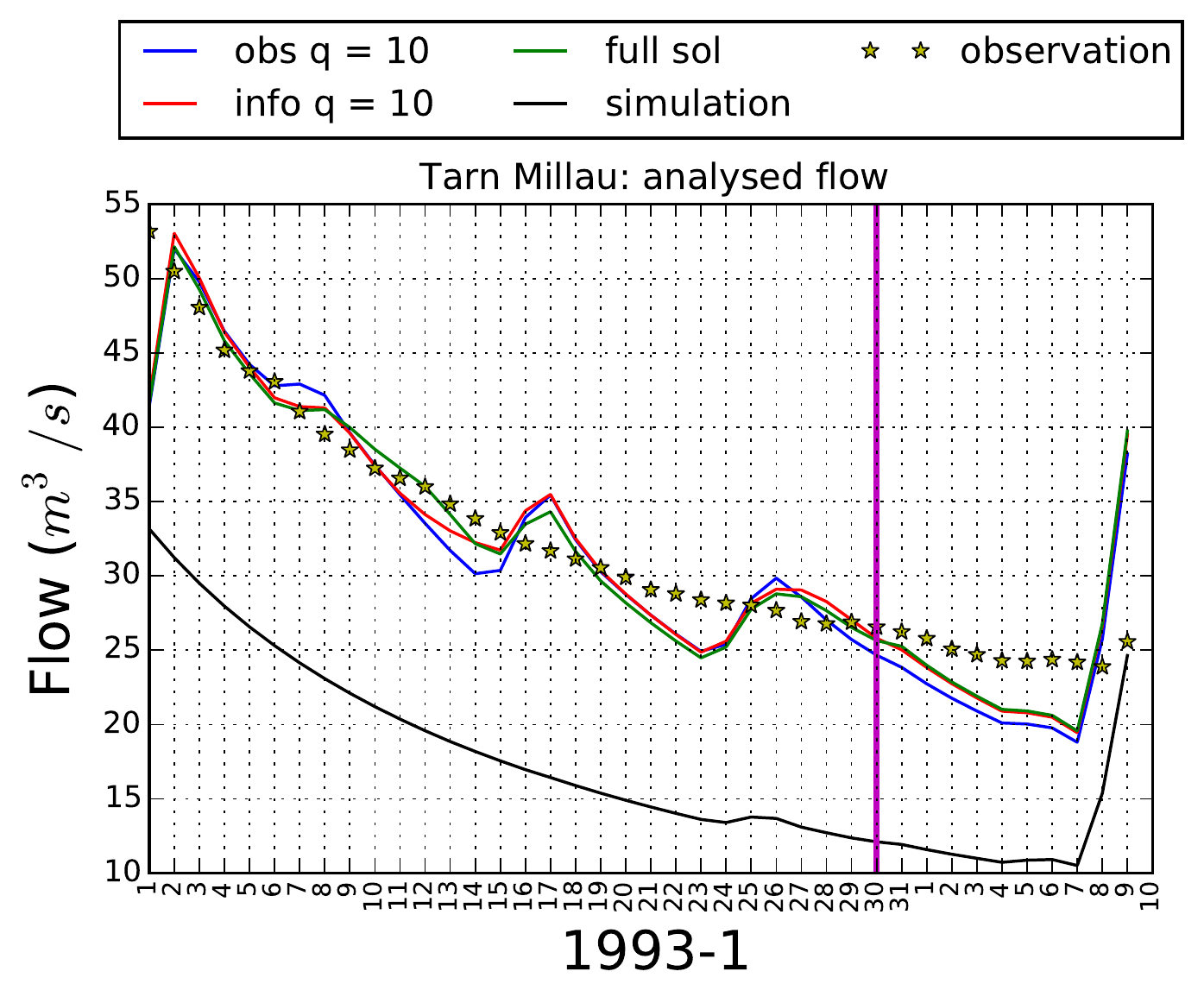}
    \includegraphics[width = 2.2  in]{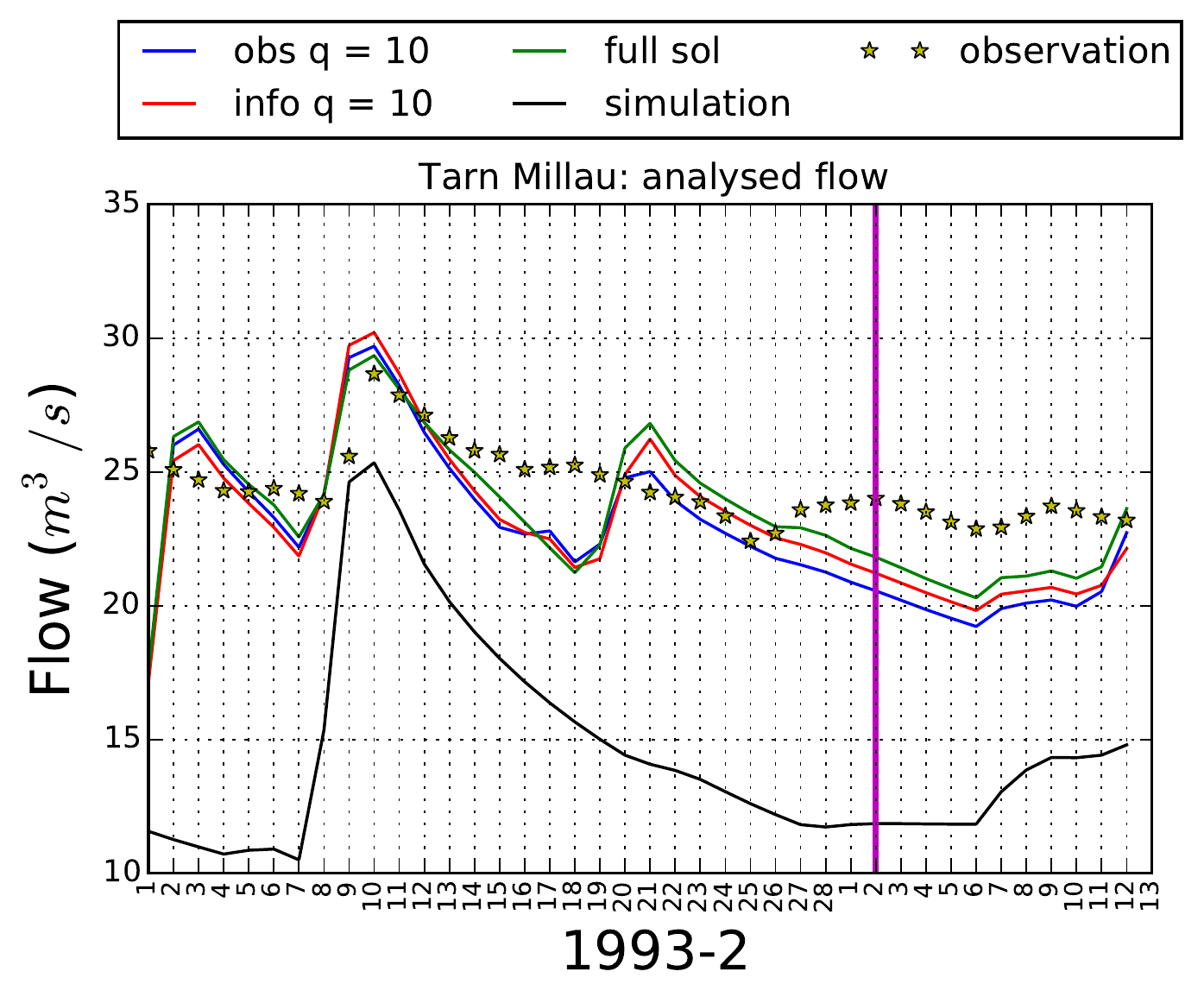}
    \includegraphics[width = 2.2  in]{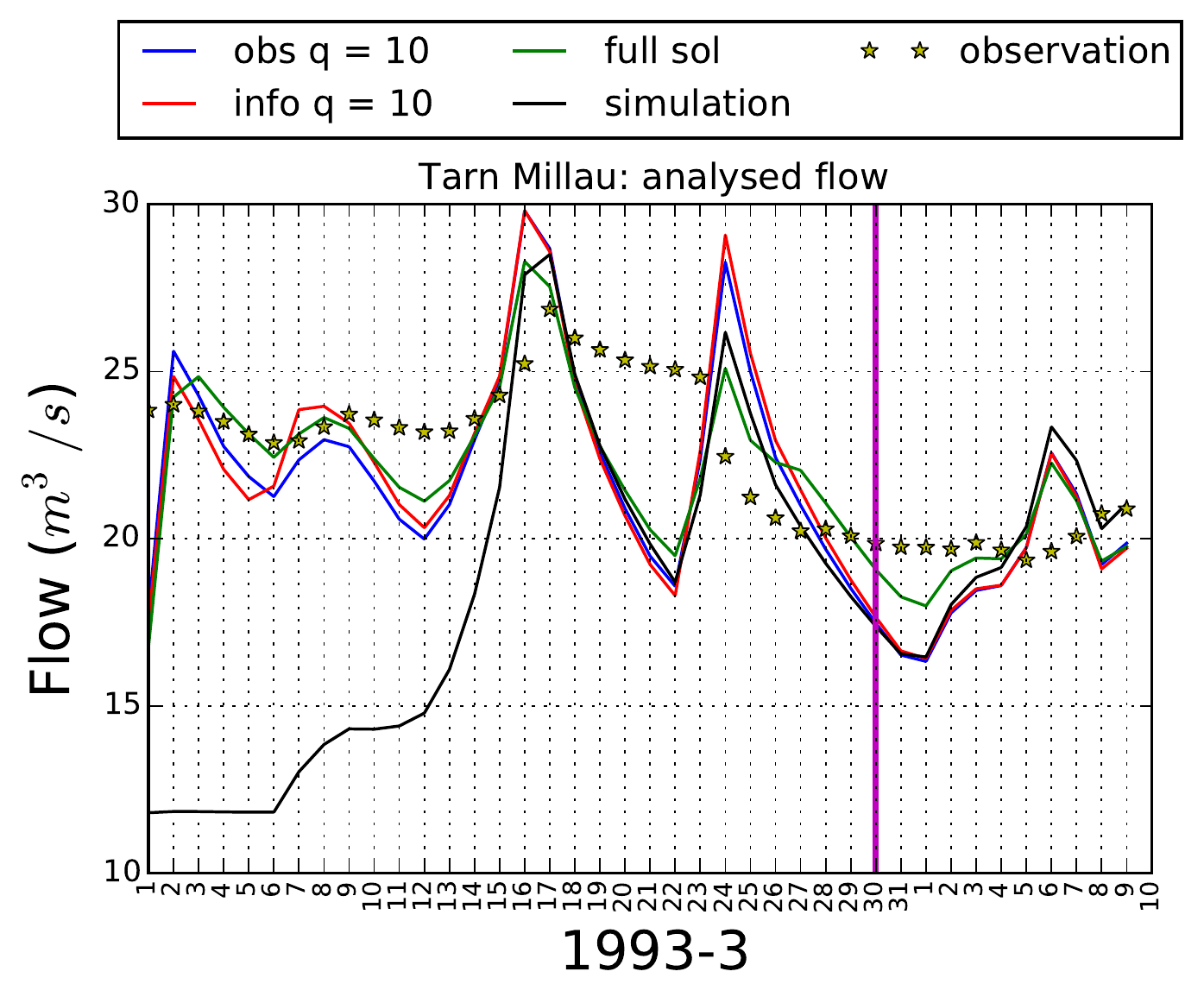}\\
    \includegraphics[width = 2.2  in]{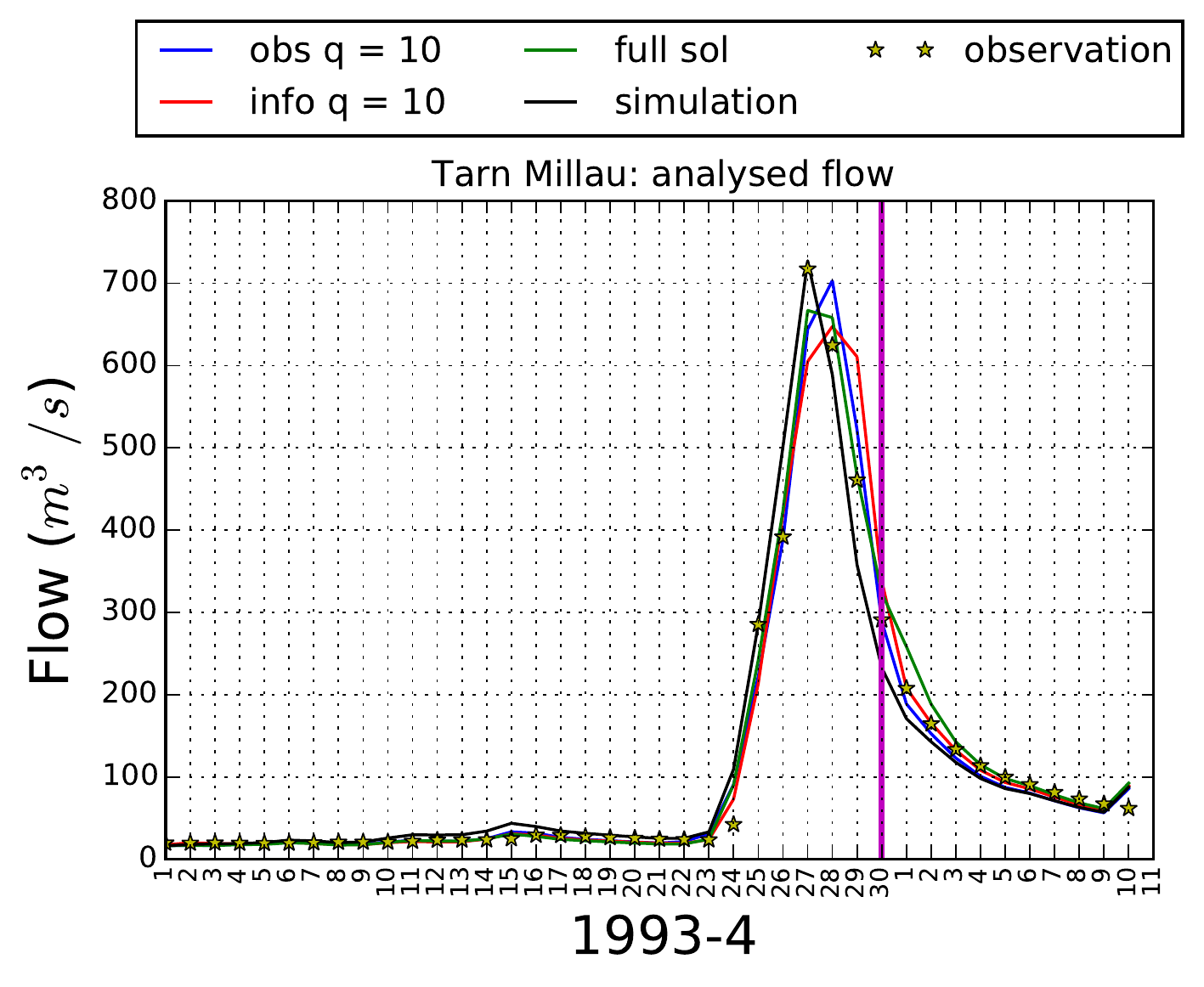}
    \includegraphics[width = 2.2  in]{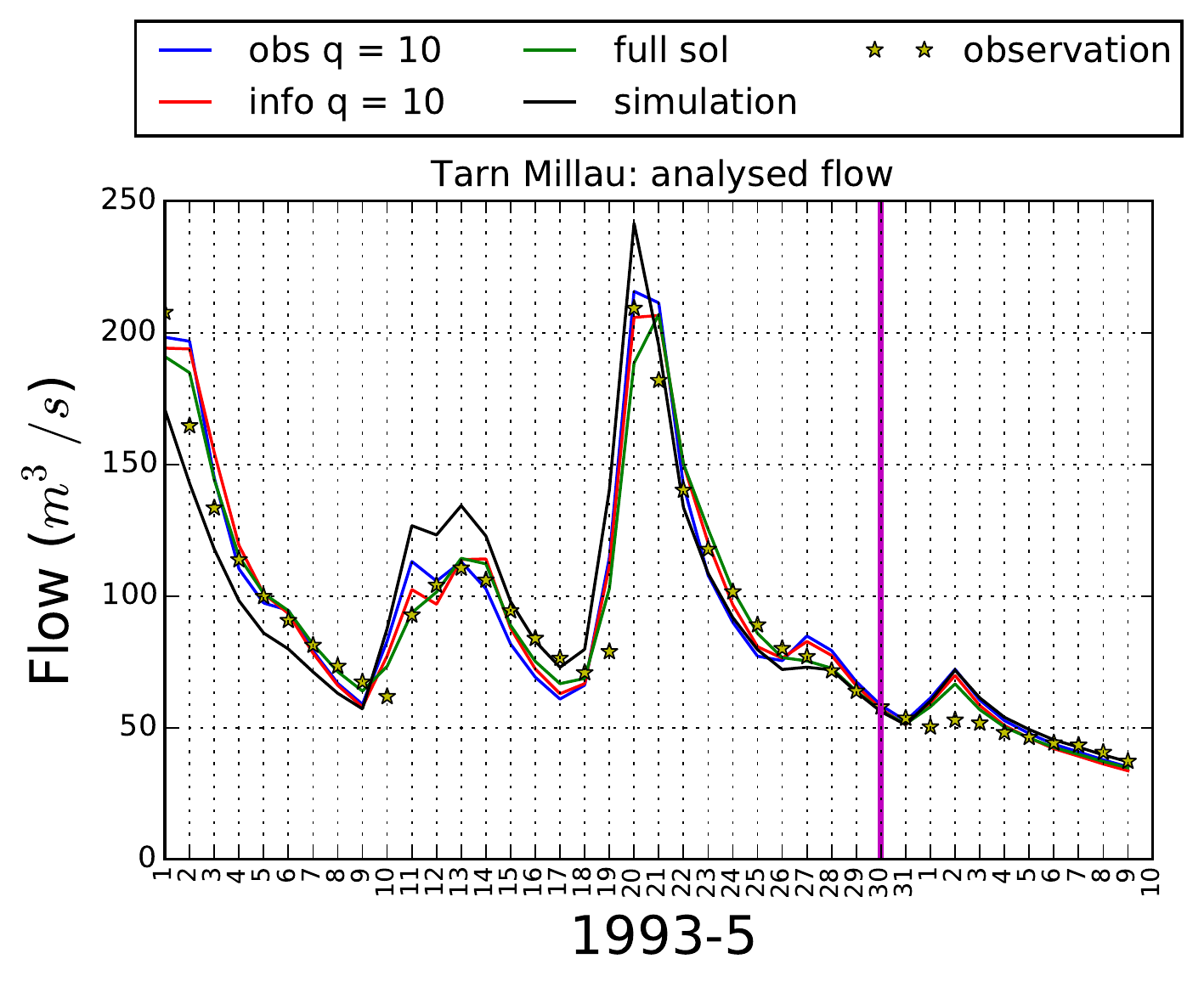}
    \includegraphics[width = 2.2  in]{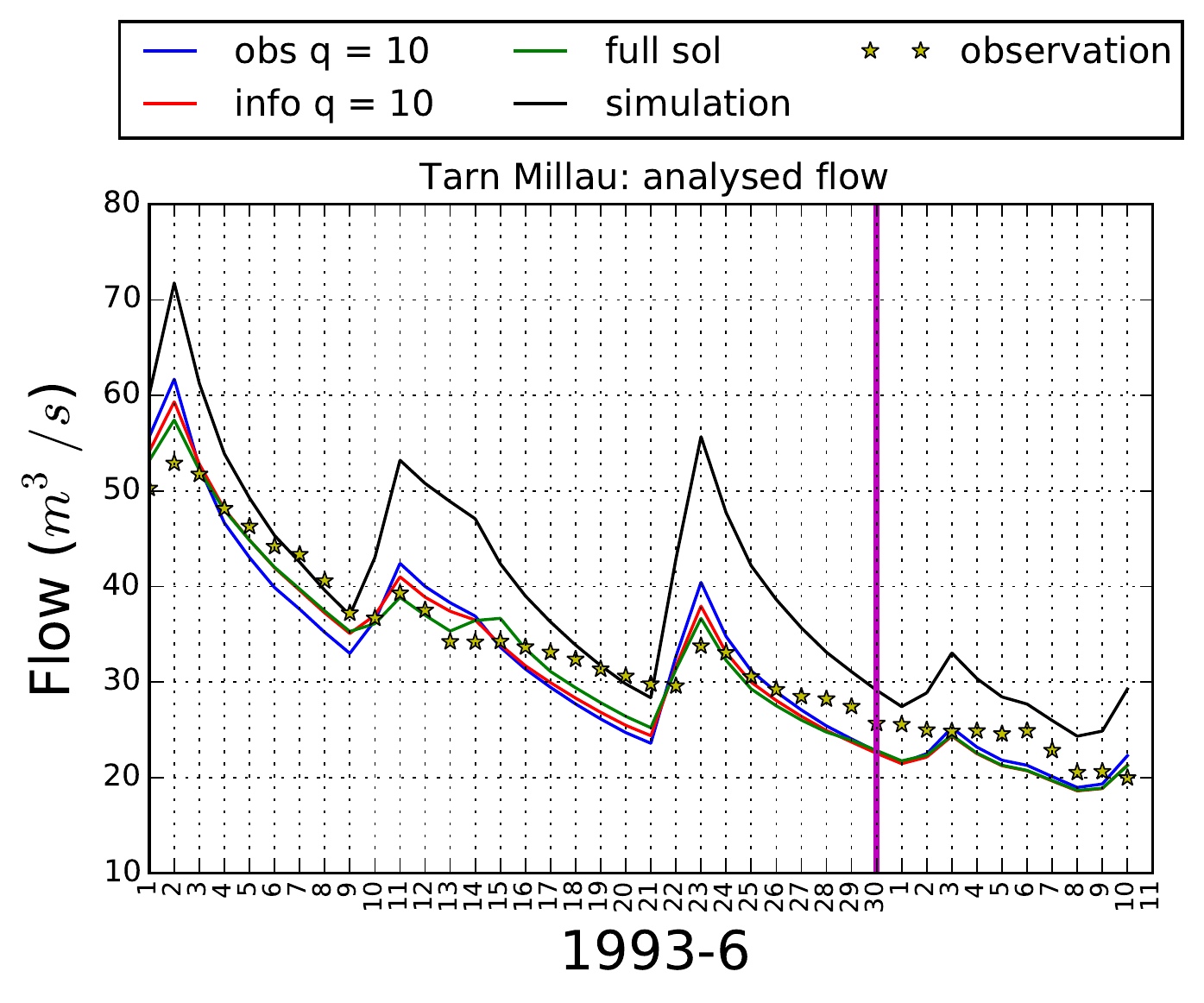}\\
    \includegraphics[width = 2.2  in]{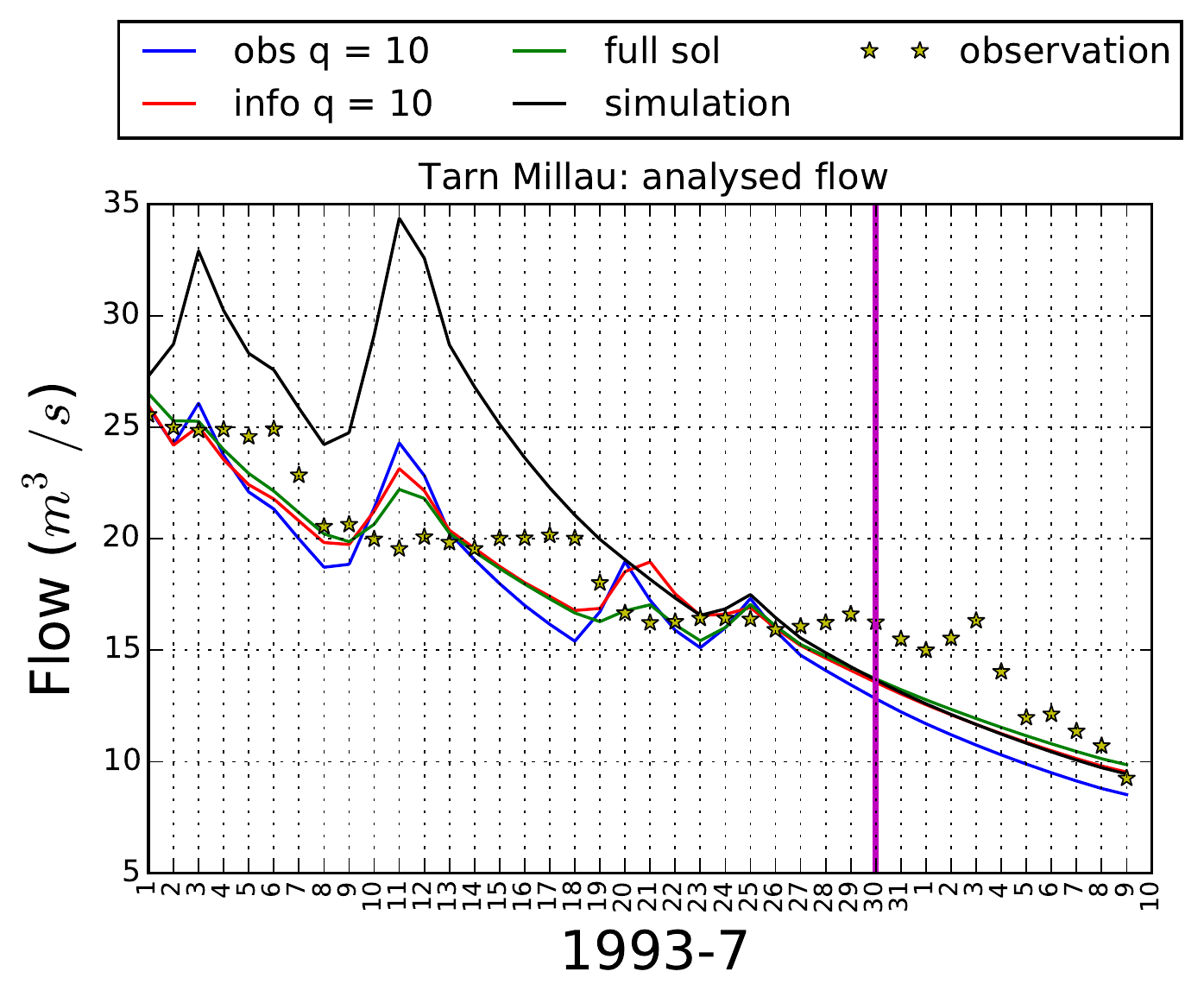}
    \includegraphics[width = 2.2  in]{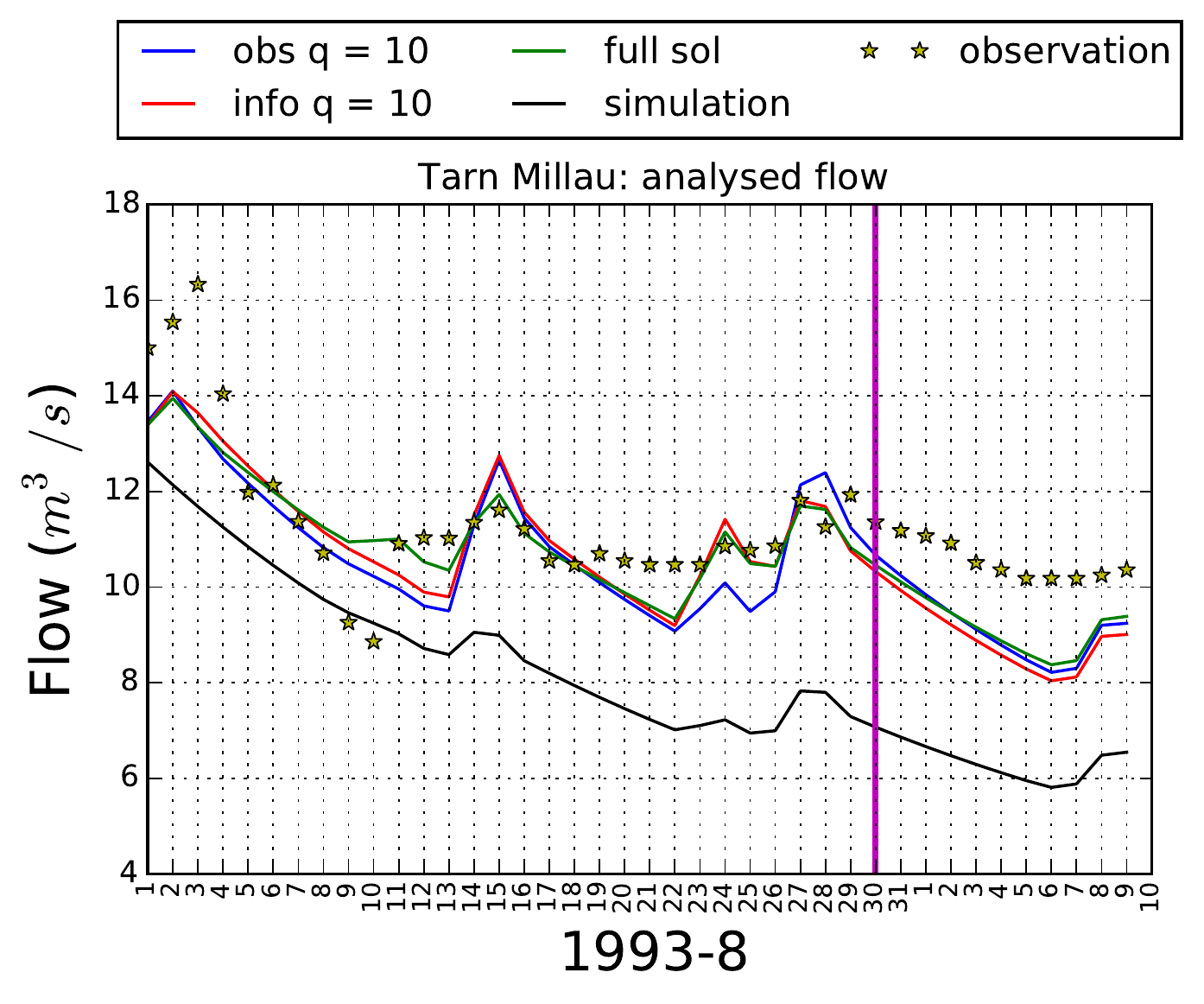}
    \includegraphics[width = 2.2  in]{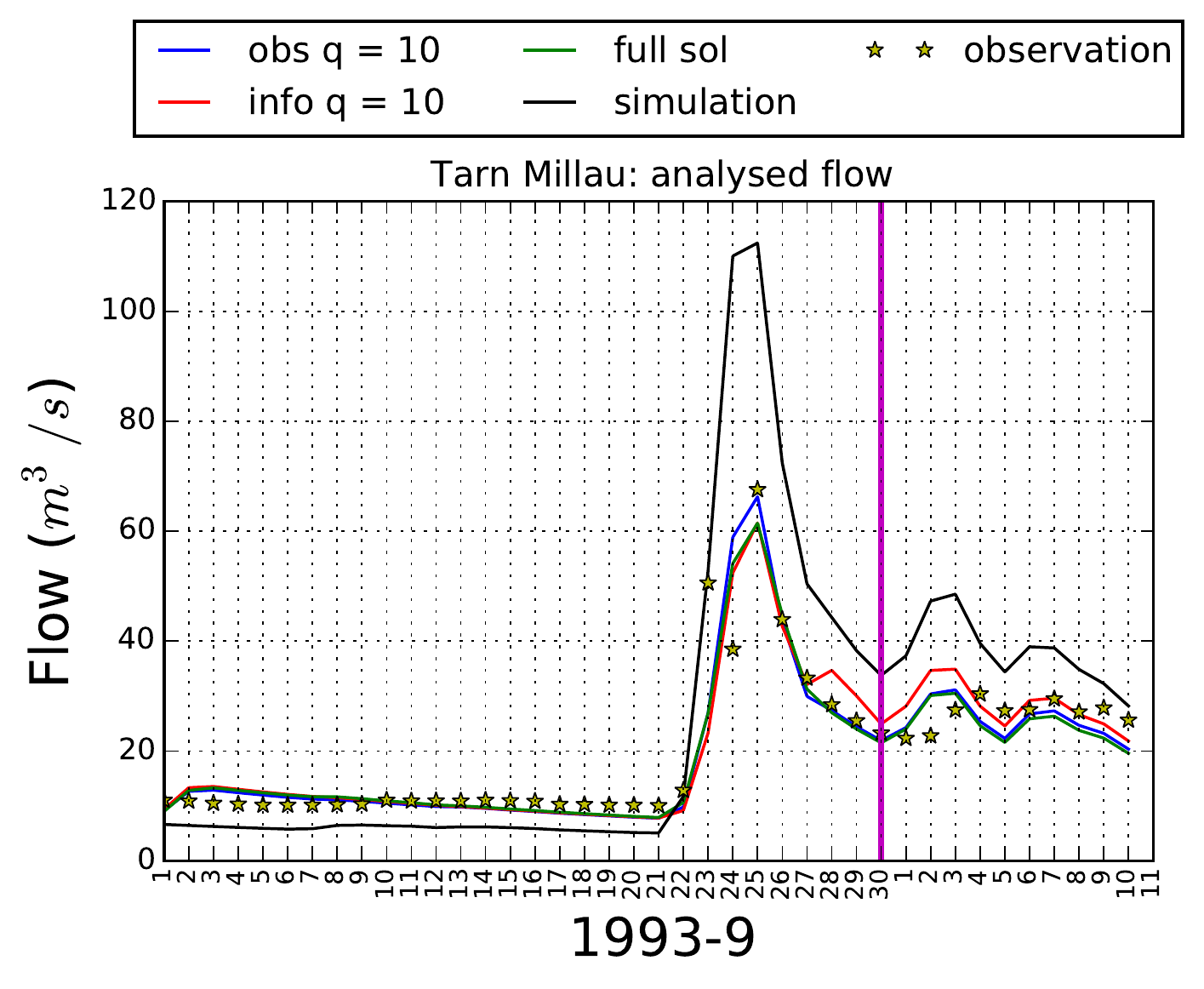}\\
    \includegraphics[width = 2.2  in]{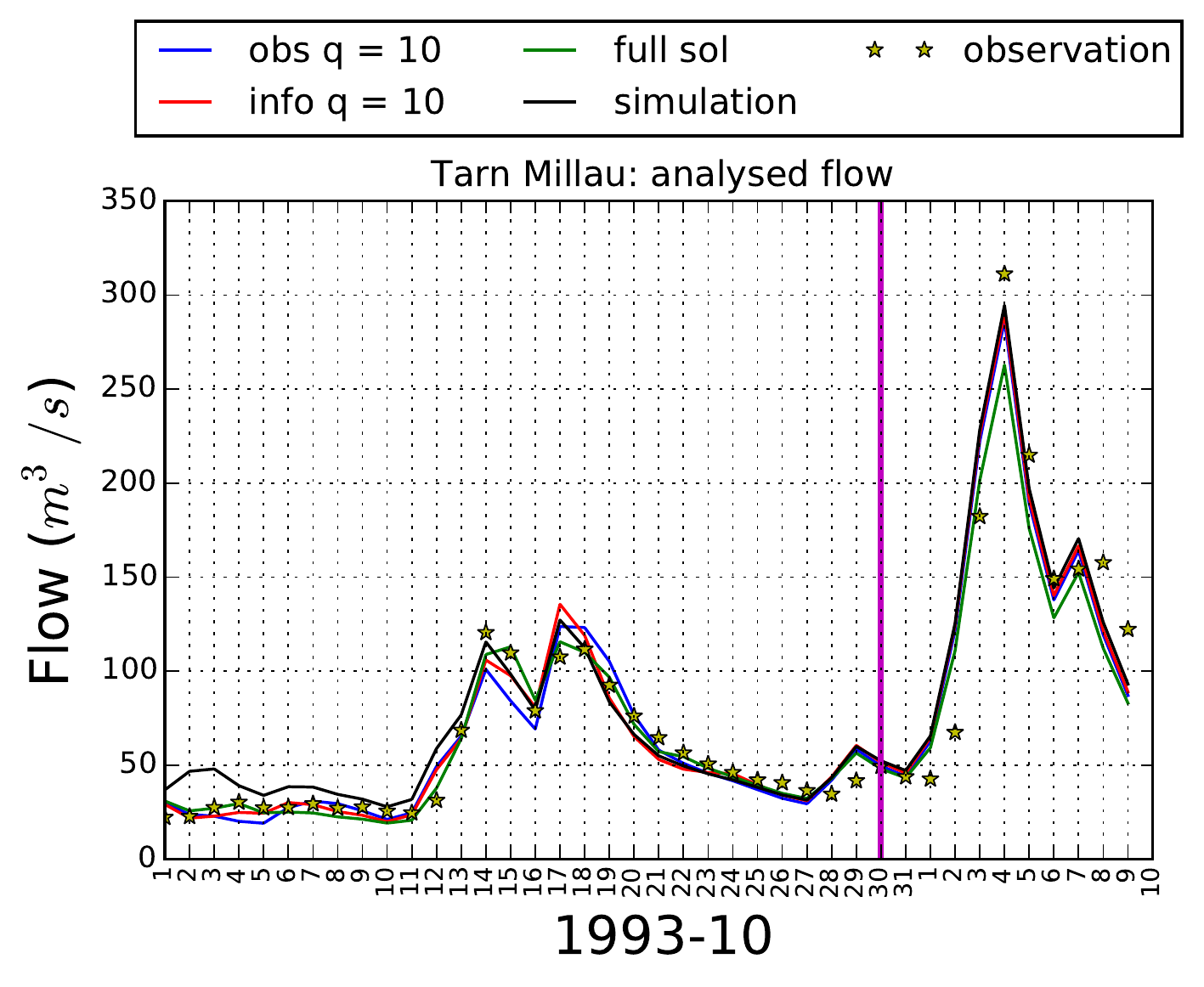}
    \includegraphics[width = 2.2  in]{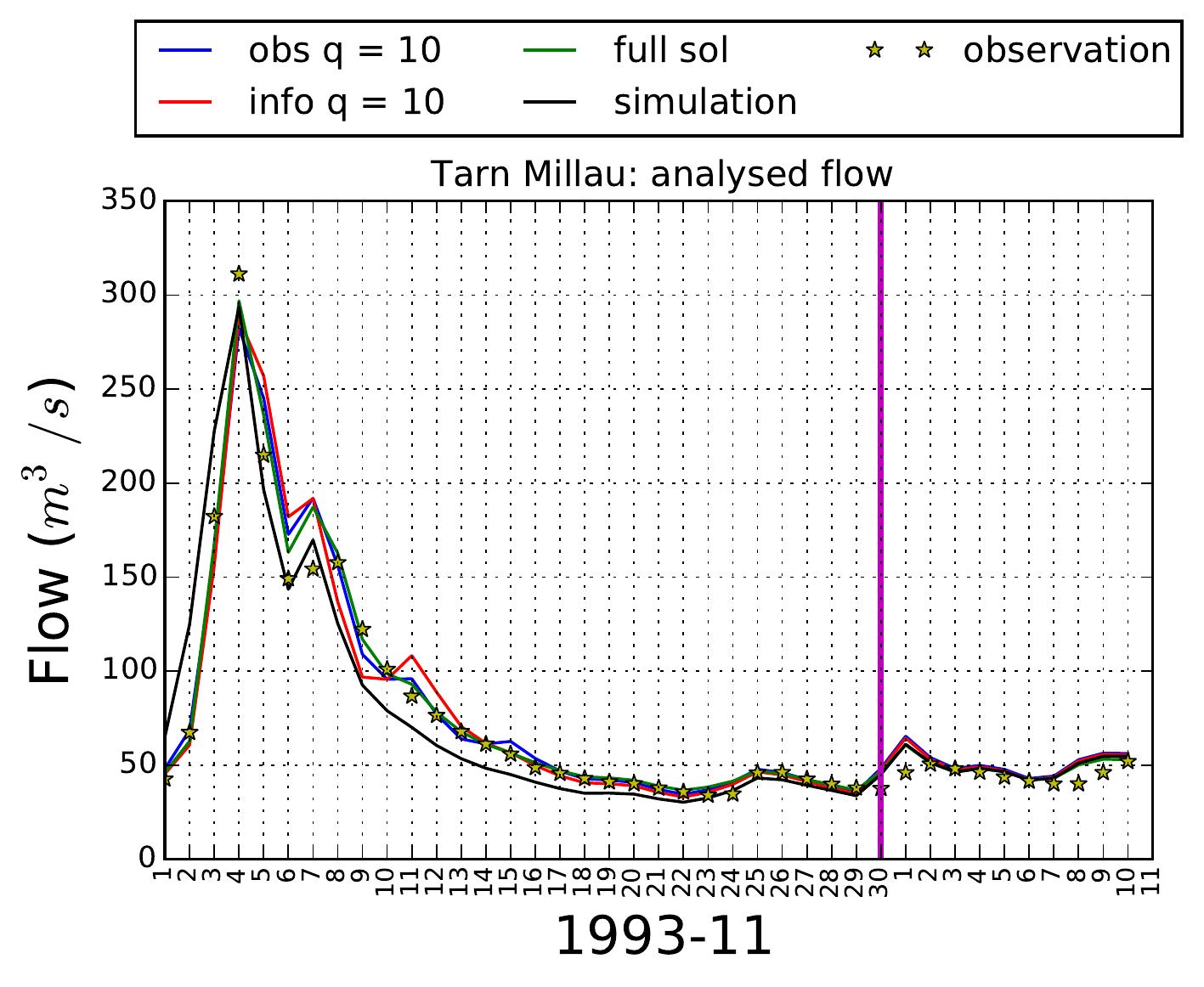}
    \includegraphics[width = 2.2  in]{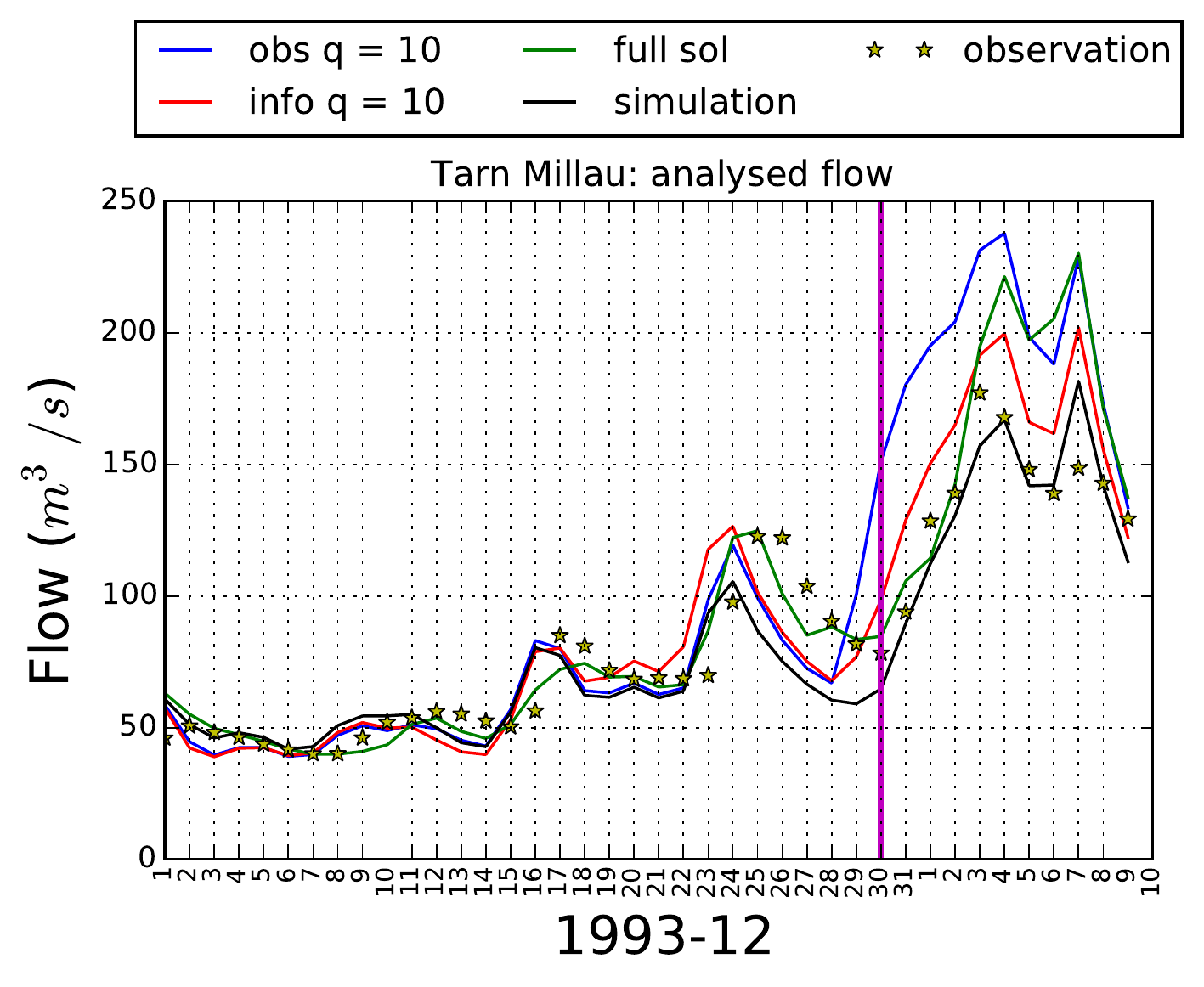}
    \caption{The reconstructed and predicted river flow at TM for OC (obs) ($q=10$), IC (info) ($q=10$) and full rank DA solutions of different months in 1993. The left side of the vertical line represents the flow reanalysis while the right side represents the prediction }
  \label{fig:evo_inno3}
\end{figure}

We draw the reconstructed river flow (i.e $\mathcal{H}(\bx_a)$) at TM of each of those 12 assimilation windows, for both corrections with compressed and full observation data in Fig. \ref{fig:evo_inno3} where the yellow stars represent the daily observations. Based on the method described in Eq.~\ref{eq:q_op}, the optimal truncation parameter reads $q_\textrm{optima l} = 22$. Here, we display the results when $q = 10$ in order to voluntary emphasize the difference between the two approaches as shown in  Fig~\ref{fig:1193-1-1 flow}.  A vertical line in each graph separates the reanalysis (left) and the prediction (right). We notice that the reconstructed curves issued from OC (blue) and IC (red) are similar in most cases, both being adequately close to the full rank assimilation (green), compared to the original simulation. Some exceptions can be found, for example, in the assimilation window of December 1993 where the prediction is covered by a flood period. It seems that the information-based approach provides a better performance, especially for flow forecasting at that moment. In general, as demonstrated in \cite{cheng2020b}, meteorological factors can impact the assimilation precision significantly. DA algorithms often perform better during drought periods (see Jun, July, August in Fig. \ref{fig:evo_inno3}) where the prior observation minus background (O-B) innovations are more consistent (i.e always being under-estimated or over-estimated). Contrarily, in flood periods where O-B innovations are usually more turbulent, more careful attention might be taken when performing compression methods.

\section{Discussion}
Sequential data assimilation algorithms can be computationally challenging, especially for large scale systems such as NWP, remote sensing, or geophysical problems.
Data compression techniques commonly used in DA problems have recently received increasing interest in reducing the computational burden. Much effort has been devoted to improving the algorithm efficiency without diminishing the accuracy of assimilation reconstruction and forecasting. Classical compression approaches consist of either extracting the principal vectors of observation dynamics or identifying the directions that contribute the most to the prior-posterior information gap. For both methods, the lack of precise knowledge on prior error covariances stands for an essential obstacle, as mentioned in several previous studies. Furthermore, the limited number of background/observation trajectories often entails a poor empirical estimation. In this paper, we have introduced a concept of observation compression benefiting from existing piecewise covariance estimation, establishing a natural connection between the posterior error covariance diagnosis and data compression techniques. More precisely, 
we assume that the error covariances (both $\bB$ and $\bR$) are flow-independent over some specific time periods, which allows an estimation based on time-variant residuals. Therefore, a much smaller number of background/observation trajectories are required for non-parametric covariance estimation. Different estimation formulations are possible depending on the prior knowledge of the $\textbf{R}$ matrix. The choice of flow-independent windows, as well as the residual sampling densities, is essential in these approaches, especially for the $\textbf{H} \textbf{B} \textbf{H}^T$ estimation. When the samplings are either too dense or too sparse, the assumptions of covariance estimation approaches might be unsatisfied, leading to a less optimal observation compression. These aspects are numerically analyzed in the twin experiments of a 2D shallow water model with non-linear dynamics with the perfect knowledge of the $\bR$. Numerical results show a significant advantage of the information-based compression in terms of assimilation accuracy, compared to the observation-based one. As for the industrial hydrological model, posterior covariance estimation which requires the knowledge of the analyzed states $\bx_a$, is needed since the $\bR$ matrix is not known \textit{a priori}. In this application, both the OC and IC compression methods rely on the flow-independent estimation of $\bR$, showing competitive performance regarding the flow reanalysis and the forecasting accuracy. A meteorological effect is also briefly discussed in this hydrological application, which indicates that different numbers of modes should be chosen in different periods of the year regarding the hydrological properties. Future work can be considered to improve the algorithm efficiency and flexibility under industrial conditions, for example, by using parametric covariance tuning methods or spatial localization techniques. Another important limitation of the current approach stands for the time invariance of the observations error covariance on some time-scale, limiting for instance the use of moving observation sensors. Indeed, if observation positions change, both the observation matrix $\mathbf{R}$ and the transformation operator $\mathcal{H}$ can not be considered flow-independent,  leading to difficulties when applying Desroziers-type methods. Future work can be considered to use interpolation approaches to construct a global observation set which includes all time-variant observation positions. Another perspective of this study could be to further examine the optimal choice of the sampling density while estimating the error covariances, for example, with the help of uncertainty quantification methods for dynamical systems.

\section*{Acknowledgement}
The authors would like to thank Dr. Bertrand Iooss and Dr. Angélique Ponçot for fruitful discussions about the compression methodology and the hydrological application. This work was supported by EDF R\&D. This research was partially funded by the Leverhulme Centre for Wildfires, Environment and Society through the Leverhulme Trust, grant number RC-2018-023. The authors are grateful to an anonymous  reviewer for the useful remarks on the manuscript.

\bibliographystyle{elsarticle-num-names}
\bibliography{main.bib}
\end{document}